\documentclass[12pt,a4paper,oneside]{article}
\usepackage{amsmath,amsfonts,amssymb,latexsym}

\newtheorem{theorem}{Theorem}[section]
\newtheorem{lemma}[theorem]{Lemma}

\newtheorem{corollary}[theorem]{Corollary}

\newtheorem{observation}[theorem]{Observation}
\newenvironment{proof}
{\bigskip\noindent{\sc Proof.}\ \ \rm }{\hfill$\Box$\bigskip}

\title{A sharp interpolation between the H\"older and Gaussian Young inequalities}

\author{Paolo Da Pelo$^1$\quad Alberto Lanconelli$^2$ \quad  Aurel I. Stan$^3$}

\date{\empty}

\begin{document}

\maketitle
\begin{center}
{\noindent
\begin{tabular}{cc}
& $^1$Dipartimento di Matematica\\
& Universita' degli Studi di Bari\\
& Via E. Orabona, 4\\
& 70125 Bari - Italia\\
& E-mail: \emph{paolo.dapelo@uniba.it}\\
\end{tabular}}
{\noindent
\begin{tabular}{cc}
& $^2$Dipartimento di Matematica\\
& Universita' degli Studi di Bari\\
& Via E. Orabona, 4\\
& 70125 Bari - Italia\\
& E-mail: \emph{alberto.lanconelli@uniba.it}\\
\end{tabular}}
{\noindent
\begin{tabular}{cc}
& $^3$Department of Mathematics\\
& Ohio State University at Marion\\
& 1465 Mount Vernon Avenue\\
& Marion, OH 43302, U.S.A.\\
& E-mail: {\em stan.7@osu.edu}\\
\end{tabular}}
\end{center}

\numberwithin{equation}{section}

\bigskip

\begin{abstract}
We prove  a very general sharp inequality of the H\"older--Young--type for functions defined on infinite dimensional Gaussian spaces.  We begin by considering a family of commutative products for functions which interpolates between the point--wise and Wick products; this family arises naturally in the context of stochastic differential equations, through Wong--Zakai--type approximation theorems, and plays a key role in some generalizations of the Beckner--type Poincar\'e inequality.  We then obtain a crucial integral representation for that family of products which is employed, together with a generalization of the classic Young inequality due to Lieb, to prove our main theorem. We stress that our main inequality contains as particular cases the H\"older inequality and Nelson's hyper-contractive estimate, thus providing a unified framework for two fundamental results of the Gaussian analysis.
\end{abstract}



\section{Introduction}

The celebrated Wong-Zakai approximation theorem \cite{WZ} establishes that  if $\{W_t^{\epsilon}\}_{t\geq 0}$ denotes a ``good" approximation of the white noise $\{W_t\}_{t\geq 0}$, then for any smooth functions $b,\sigma:\mathbb{R}\to\mathbb{R}$ the solution of the random differential equation
\begin{eqnarray}\label{WZ approximation}
\dot{X}_t^{\epsilon}=b(X_t^{\epsilon})+\sigma(X_t^{\epsilon})\cdot W_t^{\epsilon}
\end{eqnarray}
converges in probability, as $\epsilon$ goes to zero, to the solution of the stochastic differential equation
\begin{eqnarray}\label{WZ}
dX_t=b(X_t)dt+\sigma(X_t)\circ dW_t,
\end{eqnarray}
where the symbol $\circ dW_t$ denotes Stratonovich stochastic integration (recently an analogous result in the context
of stochastic partial differential equations has been obtained in the beautiful paper \cite{HP}). Replacing the point-wise product appearing
between $\sigma(X_t^{\epsilon})$ and $W_t^{\epsilon}$ in equation (\ref{WZ approximation}) with the Wick product, one gets (\cite{HO}), under the assumption of a linear diffusion
coefficient $\sigma$, the convergence, as $\epsilon$ goes to zero, to the It\^o version of equation (\ref{WZ}).\\
In the paper \cite{dpls13} the authors introduced the following family of multiplications
\begin{eqnarray}\label{tau product}
f\diamond_{\tau}g:=\tau^{-N}\Big(\tau^{N}f\cdot \tau^{N}g\Big),\quad \tau\in (0,2],
\end{eqnarray}
where $N$ denotes the number or Ornstein-Uhlenbeck operator (we refer the reader to the next section for a rigorous definition of this product in terms of second quantization operators).
This family interpolates between the point-wise product (when $\tau=1$) and the Wick product (in the limit as $\tau$ tends to zero); replacing equation (\ref{WZ approximation}) with
\begin{eqnarray*}
\dot{Y}_t^{\epsilon}=b(Y_t^{\epsilon}) + Y_t^{\epsilon}\diamond_{\tau} W_t^{\epsilon},
\end{eqnarray*}
one can prove (\cite{dpls13}) the convergence of $Y_t^{\epsilon}$ to the solution of
\begin{eqnarray*}
dY_t=b(Y_t)dt+Y_t\circ_{\tau} dW_t,
\end{eqnarray*}
where $\circ_{\tau} dW_t$ denotes stochastic integration with evaluation point at
\begin{eqnarray*}
t_i^*: =  t_{i - 1}+\frac{\tau}{2}(t_i - t_{i - 1})
\end{eqnarray*}
(this gives Stratonovich for $\tau=1$ and It\^o for $\tau=0$).  \\
The family of products defined in (\ref{tau product}) turns out to be useful also in the study of Poincar\'e-type inequalities. In fact, an important generalization of the classic
 Poincar\'e inequality (\cite{Chernoff},\cite{Nash})
\begin{eqnarray}\label{poincare}
\int f^2(w)d\mu(w)-\Big(\int f(w)d\mu(w)\Big)^2\leq\int \Vert Df(w)\Vert^2d\mu(w)
\end{eqnarray}
(here $\mu$ is a Gaussian measure defined on a possibly infinite dimensional space and $Df$ is a suitable notion of gradient of $f$) is the one proposed in \cite{Beckner} which reads for $\tau\in [0,1]$ as
\begin{eqnarray}\label{poincare2}
\int f^2(w)d\mu(w)-\int |\tau^Nf(w)|^2d\mu(w)\leq (1-\tau)\int\Vert Df(w)\Vert^2d\mu(w).
\end{eqnarray}
Inequality (\ref{poincare2}) coincides with (\ref{poincare}) for $\tau=0$ and with the logarithmic Sobolev inequality (\cite{Gross}) in the limit as $\tau$ tends to one (after an application of the Nelson's hyper-contractive estimate).
Observe that since for any $g$ one has
\begin{eqnarray*}
\int \tau^{-N}g(w)d\mu(w)=\int g(w)d\mu(w)
\end{eqnarray*}
it is possible to rewrite (\ref{poincare2}) as
\begin{eqnarray*}
\int f^2(w)d\mu(w)-\int \tau^{-N}|\tau^Nf(w)|^2d\mu(w)\leq (1-\tau)\int\Vert Df(w)\Vert^2d\mu(w)
\end{eqnarray*}
or equivalently as
\begin{eqnarray}\label{poincare tau}
\int f^2(w)d\mu(w)-\int (f\diamond_{\tau}f)(w)d\mu(w)\leq (1-\tau)\int\Vert Df(w)\Vert^2d\mu(w).
\end{eqnarray}
It has been proved in \cite{DLS14} the validity of inequality (\ref{poincare tau}) for all the probability measures obtained convolving the Gaussian measure $\mu$ with a
probability measure satisfying an exponential integrability condition.\\
It is then clear from the preceding discussion that the family of products defined in (\ref{tau product}) connects intrinsically point-wise multiplication and Stratonovich integral on one side and Wick product and It\^o integral on the other side. This connection is in addition related to the interplay between the Poincar\'e and logarithmic Sobolev inequalities. The aim of the present paper is to obtain an inequality for the
$L^r$--norm of $f\diamond_{\tau} g$ in terms of the  $L^p$--norm of $f$ and the $L^q$--norm of $g$ for suitable $p,q,r\in [1,+\infty]$ and $\tau\in [0,2]$.  We obtain a very general and sharp inequality which coincides with the classic H\"older inequality for $\tau=1$, as expected from the point of view of the interpolating nature of our family of products, and with the sharp Young--type inequality for the Wick product obtained in \cite{dpls11} for $\tau=0$. (From a probabilistic point of view, the Wick product plays in Gaussian spaces the same role played by
 the convolution in spaces equipped with the Lebesgue measure; that is why we call the inequality for the Wick product of Young--type).\\
The main purpose of this paper is to find necessary and sufficient conditions to have inequalities of the form:
\begin{eqnarray*}
\parallel f \diamond_\tau g \parallel_r & \leq & \parallel \Gamma(C)f \parallel_{\tau} \cdot \parallel \Gamma(D)g \parallel_{\tau},
\end{eqnarray*}
for all $f \in L^p(\mu)$ and $g \in L^q(\mu)$.\\
In general, as also observed in \cite{dpls11}, inequalities about the norms of Wick products are related to sharp inequalities (that means inequalities with
best constants) from classic Harmonic Analysis like: Young and Lieb inequalities, see: \cite{b1975}, \cite{bl1976}, \cite{l1990}, and \cite{ll2001}.
The sharp constant being $1$ allows us to pass to the limit as the dimension $d$ goes to infinity, having the same inequalities even in the infinite--dimensional case.\\
\par The paper is structured as follows. In section $2$ we give a minimal background on the construction of an infinite dimensional Gaussian probability space and second quantization operators.
In section 3, we review the definition of the $t$--Wick product and extend it to the definition of the $T$--Wick product, where $T$ is an
operator. We also review the definition of the exponential functions.
In section $4$ we prove an important integral representation for the Gaussian $T$--Wick products for a specific class of operators $T$.
In section 5, we use the integral representation found in section 4 and Lieb theorem from \cite{l1990}, to prove the main inequality from this paper in the dimension $d = 1$.
We use Minkowski integral inequality, to extend the inequality from dimension $d = 1$, to every finite dimension $d \geq 2$. Finally, we extend the
inequality to the infinite dimensional case.

\section{Background}
There are many ways to introduce the Gaussian Wick product and second quantization operators,
all of them being equivalent.
One can use an abstract Gel'fand triple and work with Hida White Noise Distribution Theory, see
\cite{k1996} or \cite{o1994}.
Another way is to use Malliavin Calculus, see \cite{bell06}. There is also a third way, using the
theory of Gaussian Hilbert spaces, see \cite{jan97}.
We will use Hida White Noise Distribution Theory, to make the connection with the stochastic
integral.\\
\par Let $E$ be a real separable Hilbert space, and $A$ a self--adjoint operator on $E$ having a discrete spectrum $\{\lambda_n\}_{n \geq 0}$,
such that:
\begin{enumerate}

\item There exists an orthonormal basis $\{e_n\}_{n \geq 0}$ of $E$, such that for all $n \geq 0$,
\begin{eqnarray}
Ae_n = \lambda_ne_n.
\end{eqnarray}

\item $1 < \lambda_1 < \lambda_2 < \cdots$

\item The operator $A^{-1}$ is a Hilbert--Schmidt operator.

\end{enumerate}
The inner product and norm of $E$ are denoted by $( \cdot$, $ \cdot)$ and $| \cdot |_0$, respectively.
For each $p \geq 0$, we define the norm:
\begin{eqnarray}
|f|_p^2 & := & \left|A^pf\right|_0^2\\
& = & \sum_{n = 0}^{\infty}\lambda_n^{2p}(f, e_n)^2.
\end{eqnarray}
For each $p \geq 0$, we define the space:
\begin{eqnarray}
{\mathcal E}_p & := & \{f \in E \mid |f|_p < \infty\}.
\end{eqnarray}
${\mathcal E}_p$ is a Hilbert space with norm $| \cdot |_p$.
If $0 \leq p < q$, then ${\mathcal E}_q \subset
{\mathcal E}_p$.\\
We define the space:
\begin{eqnarray}
{\mathcal E} & = & \cap_{p = 0}^{\infty}{\mathcal E}_p
\end{eqnarray}
and equip it with the locally convex topology given by the family of norms
$\{| \cdot |_p\}_{p \geq 0}$. The space ${\mathcal E}$ is a nuclear space.\\
\par For each $p \geq 0$, the dual of the space ${\mathcal E}_p$ is the space
${\mathcal E}_{-p}$,
which is the completion of the space $E$, with respect to the norm
$| \cdot |_{-p}$, defined as:
\begin{eqnarray}
|f|_{-p}^2 & := & \left|A^{-p}f\right|_0^2\\
& := & \sum_{n = 0}^{\infty}\lambda_n^{-2p}(f, e_n)^2.
\end{eqnarray}
Of course, if $0 \leq p < q$, we have:
\begin{eqnarray}
E \subset {\mathcal E}_{-p} \subset {\mathcal E}_{-q}.
\end{eqnarray}
The dual of the space ${\mathcal E}$ is the space ${\mathcal E}'$, which
can be written as:
\begin{eqnarray}
{\mathcal E}' & = & \cup_{p = 0}^{\infty}{\mathcal E}_{-p}.
\end{eqnarray}
The dual space ${\mathcal E}'$ is equipped with the inductive limit topology of the
(locally convex) topologies given by the norms $\{| \cdot |_{-p}\}_{p \geq 0}$.
We obtain in this way the following Gel'fand triple:
\begin{eqnarray}
{\mathcal E} \subset E \subset {\mathcal E}'.
\end{eqnarray}
By Minlos theorem there exists a unique probability measure $\mu$ on the dual space
${\mathcal E}'$ of ${\mathcal E}$, such that, for all
$\xi \in {\mathcal E}$, we have:
\begin{eqnarray}
\int_{{\mathcal E}'}e^{i\langle x, \xi \rangle}d\mu(x) & = & e^{-(1/2)|\xi|_0^2}, \label{minlos}
\end{eqnarray}
where $\langle \cdot$, $\cdot \rangle$ denotes the bilinear pairing of ${\mathcal E}'$ and
${\mathcal E}$, see page 16 of \cite{k1996}.
Formula (\ref{minlos}) says that as a random variable the continuous function $\langle \cdot$,
$\xi \rangle$ is normally distributed with mean $0$ and variance $|\xi|_0^2$, for every
$\xi \in {\mathcal E}$. This observation is very important, since by approximating in the norm $| \cdot |_0$ of $E$,
every element $f$ of $E$, by a sequence $\{\xi_n\}_{n \geq 1}$ of elements of ${\mathcal E}$, we obtain
a Cauchy sequence $\{\langle \cdot$, $\xi_n \rangle\}_{n \geq 1}$ in $L^2({\mathcal E}'$, $\mu)$ of normally
distributed random variables. The $L^2$--limit of this sequence is denoted by $\langle \cdot$, $f \rangle$ and is
a normally distributed random variable with mean $0$ and variance $|f|_0^2$.\\
\par For every real Hilbert space $H$, we denote by $H_c$ its complexification.
We define the {\em trace operator} $\tau$ as the following element of
$({\mathcal E}_c')^{\hat{\otimes}2}$, where $\hat{\otimes}$ denotes the symmetric tensor product:
\begin{eqnarray}
\langle \tau, \xi \otimes \eta \rangle & := & \langle \xi, \eta \rangle,
\end{eqnarray}
for all $\xi$ and $\eta$ in ${\mathcal E}_c$.
We define the {\em Wick tensor} $:x^{\otimes n}:$, for every $x \in {\mathcal E}'$ as:
\begin{eqnarray*}
:x^{\otimes n}: & = & \sum_{k = 0}^{[n/2]}{n \choose {2k}}(2k - 1)!!(-1)^kx^{\otimes(n - 2k)} \hat{\otimes}\tau^{\otimes k}.
\end{eqnarray*}
If we denote by $(L^2)$ the space of all complex valued square integrable functions defined on $({\mathcal E}'$, $\mu)$, then for every function $\varphi$ in $(L^2)$,
there exists a unique sequence $\{f_n\}_{n \geq 0}$, where for all $n \geq 0$, $f_n \in E_c^{\hat{\otimes} n}$, such that:
\begin{eqnarray}
\varphi(x) & = & \sum_{n = 0}^{\infty}\langle :x^{\otimes n}:, f_n\rangle. \label{Wiener_Ito_decomposition}
\end{eqnarray}
Moreover, the square of the $(L^2)$--norm of $\varphi$ is:
\begin{eqnarray}
\parallel \varphi \parallel_0^2 & = & \sum_{n = 0}^{\infty}n!|f_n|_0^2\\
& < & \infty \nonumber,
\end{eqnarray}
where $|f_n|_0$ denotes the norm of $f_n$ computed in the space $E_c^{\otimes n}$.\\
If $B$ is a densely defined operator on $E$, and $\varphi$ is given by (\ref{Wiener_Ito_decomposition}), then
we define the {\em second quantization operator} of $B$, as:
\begin{eqnarray}
\Gamma(B)\varphi(x) & := & \sum_{n = 0}^{\infty}\langle :x^{\otimes n}:, B^{\otimes n}f_n \rangle.
\end{eqnarray}
It is not hard to see that if $B$ is a bounded operator on $E$, of operatorial norm $\parallel B \parallel \leq 1$,
then $\Gamma(B)$ is a bounded operator on $(L^2)$ of operatorial norm $\parallel \Gamma(B) \parallel = 1$.\\
In particular, if we take $B := A$, the unbounded operator used to define the Gel'fand triple
${\mathcal E} \subset E \subset {\mathcal E}'$, then
the second quantization operator $\Gamma(A)$ has properties similar to those of the operator $A$:
\begin{itemize}
\item $\Gamma(A)$ has positive eigenvalues and a set of eigenfunctions that forms an orthogonal
basis of $(L^2)$.
\item $\Gamma(A)^{-1}$ is a bounded operator.
\item For every $p > 1$, $\Gamma(A)^{-p}$ is a Hilbert--Schmidt operator on $(L^2)$.
\end{itemize}
Repeating the same constructions as before, with $(L^2)$ and $\Gamma(A)$ replacing $E$ and $A$, respectively,
we obtain a new Gel'fand triple:
\begin{eqnarray}
({\mathcal E}) \subset (L^2) \subset ({\mathcal E})^*.
\end{eqnarray}
$({\mathcal E})$ is called the space of {\em test functions}, while $({\mathcal E})^*$ is named the space of
{\em generalized functions} (or {\em Hida distributions}).\\
The bilinear pairing between $({\mathcal E})^*$ and $({\mathcal E})$ is denoted by $\langle\langle \cdot$,
$\cdot \rangle\rangle$. It must be mentioned that while the spaces involved in the first Gel'fand triple:
\begin{eqnarray*}
{\mathcal E} \subset E \subset {\mathcal E}'
\end{eqnarray*}
are vector spaces over ${\mathbb R}$, the spaces used in the second Gel'fand triple:
\begin{eqnarray*}
\left({\mathcal E}\right) \subset \left(L^2\right) \subset \left({\mathcal E}\right)^*
\end{eqnarray*}
are vector spaces over ${\mathbb C}$.\\
\par The following two theorems can be found in \cite{o1994}, pages 35--36.
\begin{theorem}
Let $\phi \in (L^2)$ have the following Wiener--It\^o expansion:
\begin{eqnarray*}
\phi(x) & = & \sum_{n = 0}^{\infty}\langle :x^{\otimes n}:, f_n \rangle,
\end{eqnarray*}
where $x \in {\mathcal E}'$, and for each $n \geq 0$, $f_n \in E_c^{\hat{\otimes} n}$, such that:
\begin{eqnarray*}
\sum_{n = 0}^{\infty}n!|f_n|_0^2 & < & \infty.
\end{eqnarray*}
Then $\phi \in ({\mathcal E})$ if and only if, for all $n \geq 0$, $f_n \in {\mathcal E}_c^{\hat{\otimes} n}$, and
for all $p \geq 0$:
\begin{eqnarray}
\parallel \phi \parallel_p^2 & := & \sum_{n = 0}^{\infty}n!|f_n|_p^2\\
& < & \infty.
\end{eqnarray}
\end{theorem}
\begin{theorem}
For each $\phi \in ({\mathcal E})^*$ there exists a unique sequence $\{F_n\}_{n \geq 0}$,
such that, for all $n \geq 0$, $F_n \in {{\mathcal E}'_c}^{\hat{\otimes} n}$:
\begin{eqnarray}
\langle\langle \phi, \varphi \rangle\rangle & = & \sum_{n = 0}^{\infty}n!\langle F_n, f_n \rangle,
\label{gendef}
\end{eqnarray}
for all $\varphi \in ({\mathcal E})$, where $\varphi$ and $\{f_n\}_{n = 0}^{\infty}$ are related by the previous theorem.\\
Conversely, given a sequence $\{F_n\}_{n = 0}^{\infty}$, such that, for each $n \geq 0$,
$F_n \in {{\mathcal E}'_c}^{\hat{\otimes} n}$ and there exists $p \geq 0$, such that:
\begin{eqnarray}
\sum_{n = 0}^{\infty}n!|F_n|_{-p}^2 & < & \infty,
\end{eqnarray}
a generalized functional $\phi \in (E)^*$ is defined by (\ref{gendef}).
In this case, we write:
\begin{eqnarray}
\phi(x) & := & \sum_{n = 1}^{\infty}\langle :x^{\otimes n}:, F_n \rangle.
\end{eqnarray}
\end{theorem}
If $\phi(x) = \sum_{n = 0}^{\infty}\langle :x^{\otimes n}:$, $F_n \rangle$, and there
exists $N \in {\mathbb N}$, such that, for all $n \geq N$, $F_n = 0$, then we call
$\phi$ a {\em polynomial} generalized function.\\
\par As a particular example of the above general construction we present the following.
Let $k$ be a natural number and $E := L^2({\mathbb R}^k$, $dx)$, where $dx$ denotes the Lebesgue
measure. We consider the following self--adjoint operator on $E$:
\begin{eqnarray}
A := \left(x_1^{2} - \frac{d^2}{dx_1^2} + 1\right)\left(x_2^{2} - \frac{d^2}{dx_2^2} + 1\right)
\cdots \left(x_k^{2} - \frac{d^2}{dx_k^2} + 1\right).
\end{eqnarray}
Then $A$ satisfies the conditions required by our construction. In this case, the nuclear space
${\mathcal E}$ becomes the Schwartz space of rapidly decreasing smooth functions, and its dual
${\mathcal E}'$ the space of tempered distributions.\\
Let ${\mathcal B}_f$ denote the set of all Borel subsets of ${\mathbb R}^k$
of finite Lebesgue measure.\\
Since for every set $X$ in ${\mathcal B}_f$, its characteristic function
$1_X$ belongs to $L^2({\mathbb R}^k$, $dx)$, we can define the $(L^2)$ random variable:
\begin{eqnarray}
B_X & := & \langle \cdot, 1_X \rangle.
\end{eqnarray}
Then the family of random variables $\{B_X\}_{X \in {\mathcal B}_f}$ is a Brownian sheet.\\
In particular, if $k = 1$, and for every $t \geq 0$, we define:
\begin{eqnarray}
B_t & := & \langle \cdot, 1_{[0, t]} \rangle,
\end{eqnarray}
then $\{B_t\}_{t \geq 0}$ is a Brownian motion process.\\
The derivative of the Brownian motion is the following polynomial generalized
function:
\begin{eqnarray}
\dot{B}_t & = & \langle \cdot, \delta_t \rangle,
\end{eqnarray}
where $\delta_t$ denotes the Dirac delta measure, for all $t \in {\mathbb R}$.

\section{Generalized Wick products and Exponential Functions}

For any non--negative integer $n$, let us denote by ${\mathcal G}_n$, the following closed subspace
of $(L^2)$:
\begin{eqnarray}
{\mathcal G}_n & := & \{\langle :x^{\otimes n}:, f_n \rangle \mid f_n \in E_c^{\hat{\otimes} n}\}.
\end{eqnarray}
We call ${\mathcal G}_n$ the space of {\em homogenous polynomial random variables of degree} $n$.
It is clear that the spaces $\{{\mathcal G}_n\}_{n \geq 0}$ are mutually orthogonal.\\
For all $n \geq 0$, let us define:
\begin{eqnarray}
{\mathcal F}_n & := & \sum_{k = 0}^n{\mathcal G}_k.
\end{eqnarray}
Let $P_n$ and $P_{< n}$ denote the orthogonal projection of $(L^2)$ onto
the closed subspaces ${\mathcal G}_n$ and ${\mathcal F}_{n - 1}$, respectively.
\\
If $\varphi(x) = \sum_{n = 0}^{\infty}\langle :x^{\otimes n}:$, $ F_n\rangle$ and
$\psi(x) = \sum_{n = 0}^{\infty}\langle :x^{\otimes n}:$, $ G_n\rangle$, with $F_n$ and $G_n$ in ${\mathcal E}_c'^{\hat{\otimes} n}$
for all $n \geq 0$, then we define the {\em classic Wick product}, $\varphi \diamond \psi$, of $\varphi$ and $\psi$, as the following
generalized function:
\begin{eqnarray}
(\varphi \diamond \psi)(x) & := & \sum_{k = 0}^{\infty}\langle :x^{\otimes k}:, h_k\rangle,
\end{eqnarray}
where
\begin{eqnarray}
h_k & := & \sum_{u + v = k}F_u \hat{\otimes} G_v.
\end{eqnarray}
It is shown in \cite{k1996}, Theorem 8.12, page 92, that the Wick product is continuous from $({\mathcal E})^* \times ({\mathcal E})^*$
into $({\mathcal E})^*$.\\
\par For every $t > 0$ (and later on we will restrict $t$ to $(0$, $2]$), the
{\em generalized Wick product} or $t$--{\em Wick product},
introduced by Da Pelo and Lanconelli (see \cite{dpl12}), can be defined
using the second quantization operator of $\sqrt{t}$ times the identity
operator $I$ as:
\begin{eqnarray}
\varphi \diamond_t \psi & := & \Gamma\left(\frac{1}{\sqrt{t}}I\right)
\left[\Gamma(\sqrt{t}I)\varphi \cdot \Gamma(\sqrt{t}I)\psi\right], \label{t_Wick_def}
\end{eqnarray}
for every $(\varphi$, $\psi)$ in a dense subspace $V_t$ of
$(L^2) \times (L^2)$,
for which\\ $\Gamma((1/\sqrt{t})I)
[\Gamma(\sqrt{t}I)\varphi \cdot \Gamma(\sqrt{t}I)\psi]$ belongs to
$(L^2)$. Such a space $V_t$ can be taken as the vector space spanned by exponential functions 
(which will be defined later).\\
We know for sure that for every two polynomial random variables $\varphi$ and $\psi$
in $(L^2)$, $\varphi \diamond_t \psi$ is also a polynomial random variable. Moreover,
$\varphi \diamond_t \psi$ can be viewed as a polynomial in the variable $t$ with
coefficients in the spaces $E_c^{\hat{\otimes} n}$, for $n \geq 0$, whose constant term
(i.e., the term without $t$) is the classic Wick product $\varphi \diamond \psi$.
To understand this, let us write:
$\varphi = \sum_{p = 0}^m f_p$
and
$\psi = \sum_{q = 0}^n g_q$,
where for all $0 \leq p \leq m$, $f_p \in {\mathcal G}_p$, and for all $0 \leq q \leq n$, $g_q \in
{\mathcal G}_q$.
Since, for all $(p$, $q) \in \{0$, $1$, $\dots$, $m\} \times \{0$, $1$, $\dots$, $n\}$,
$f_p \cdot g_q \in {\mathcal F}_{p + q}$, and the Gaussian probability measure is symmetric, we have:
\begin{eqnarray}
f_p \cdot g_q & = & P_{p + q}\left(f_p \cdot g_q\right) + P_{< (p + q)}\left(f_p \cdot g_q\right)\\
& = & f_p \diamond g_q + \sum_{\substack{k < p + q\\ k \equiv (p + q)({\rm mod} \ 2)}}P_k\left(f_p \cdot g_q\right),
\end{eqnarray}
Thus, for all $t \in (0$, $2]$, we have:
\begin{eqnarray}
\varphi \diamond_t \psi & = & \Gamma\left(\frac{1}{\sqrt{t}}I\right)
\left[\Gamma(\sqrt{t}I)\varphi \cdot \Gamma(\sqrt{t}I)\psi\right] \nonumber\\
& = & \Gamma\left(\frac{1}{\sqrt{t}}I\right)\left[\sum_{p = 0}^m t^{p/2}f_p \cdot
\sum_{q = 0}^n t^{q/2}g_q\right] \nonumber\\
& = & \Gamma\left(\frac{1}{\sqrt{t}}I\right)\left[\sum_{k = 0}^{m + n}t^{k/2}\sum_{p + q = k}f_p \cdot g_q\right]\nonumber\\
& = & \sum_{k = 0}^{m + n}t^{k/2}\sum_{p + q = k}\Gamma\left(\frac{1}{\sqrt{t}}I\right)\left(f_p \diamond g_q +
\sum_{\substack{r < k\\r \equiv k({\rm mod} \ 2)}}P_r\left(f_p \cdot g_q\right)\right). \nonumber
\end{eqnarray}
Since for all $l \geq 0$ and $h \in {\mathcal G}_l$, we have $\Gamma(1/(\sqrt{t})I)h = t^{-l/2}h$, we obtain:
\begin{eqnarray}
\varphi \diamond_t \psi & = & \sum_{k = 0}^{m + n}t^{k/2}\sum_{p + q = k}\left[t^{-k/2}f_p \diamond g_q +
\sum_{\substack{r < k\\r \equiv k({\rm mod} \ 2)}}t^{-r/2}P_r\left(f_p \cdot g_q\right)\right] \nonumber\\
& = & \sum_{k = 0}^{m + n}\sum_{p + q = k}f_p \diamond g_q + \sum_{k = 1}^{m + n}\sum_{p + q = k}
\sum_{\substack{r < k\\r \equiv k({\rm mod} \ 2)}}t^{(k - r)/2}P_r\left(f_p \cdot g_q\right) \nonumber\\
& = & \varphi \diamond g + \sum_{k = 1}^{m + n}\sum_{p + q = k}
\sum_{\substack{r < k\\r \equiv k({\rm mod} \ 2)}}t^{(k - r)/2}P_r\left(f_p \cdot g_q\right).
\nonumber
\end{eqnarray}
It follows from the last relation that, at least in the case of polynomial random variables,
the classic Wick product can be understood as the $0$--Wick product in the sense of Da Pelo and Lanconelli.
That means:
\begin{eqnarray}
\varphi \diamond \psi & = & \varphi \diamond_0 \psi\\
& := & \lim_{t \to 0^+}\varphi \diamond_t \psi.
\end{eqnarray}
For this reason, from now on, we will take $t \geq 0$, when speaking about the family of $t$--Wick products.\\
\par We can generalize this product in the following way.
For every bounded self--adjoint operator $T$ on $E$, that commutes with the operator $A$ used to define the Gel'fand triple
${\mathcal E} \subset E \subset {\mathcal E}'$,
such that $T > 0$, we
define the $T$--Wick product as:
\begin{eqnarray}
\varphi \diamond_T \psi & := & \Gamma\left(\frac{1}{\sqrt{T}}\right)
\left[\Gamma(\sqrt{T})\varphi \cdot \Gamma(\sqrt{T})\psi\right], \label{t_Wick_def}
\end{eqnarray}
for every $(\varphi$, $\psi)$ in a dense subspace $V_T$ of
$(L^2) \times (L^2)$,
for which\\
$\Gamma((1/\sqrt{T}))[\Gamma(\sqrt{T})\varphi \cdot \Gamma(\sqrt{T})\psi]$ belongs to
$(L^2)$. Such a space $V_T$ can be taken as the vector space spanned by exponential functions.\\
Here when we speak of $\sqrt{T}$ and $1/\sqrt{T}$, we understand the self--adjoint operators whose
eigenvectors are $\{e_n\}_{n \geq 0}$ (the same as the eigenvectors of $A$) and whose
eigenvalues are $\{\sqrt{\mu_n}\}_{n \geq 0}$ and $\{1/\sqrt{\mu_n}\}_{n \geq 0}$, respectively, where
$\{\mu_n\}_{n \geq 0}$ are the eigenvalues of the operator $T$.
It is clear that the $T$--Wick product is commutative and associative.\\
As before, if $T \geq 0$, by gently passing to the limit as $\epsilon \to 0^+$, we can
define the $T$--Wick product as:
\begin{eqnarray}
\varphi \diamond_T \psi & := & \lim_{\epsilon \to 0^+}\varphi \diamond_{T_{\epsilon}} \psi,
\end{eqnarray}
where
\begin{eqnarray}
T_{\epsilon} & := & T + \epsilon P_{Ker(T)},
\end{eqnarray}
where
\begin{eqnarray}
Ker(T) & := & \{x \in E \mid Tx = 0_E\}
\end{eqnarray}
and $P_{Ker(T)}$ denotes the orthogonal projection of $E$ onto $Ker(T)$.\\
\par We recall now an important family of random variables called the {\em (renormalized) exponential functions}.\\
For every $\xi \in E_c$, we define the {\em exponential
function} $\varphi_{\xi}$ {\em generated by} $\xi$ as:
\begin{eqnarray}
\varphi_{\xi} & := & \sum_{n = 0}^{\infty}\frac{1}{n!}
\langle \cdot , \xi \rangle^{\diamond n}, \label{expdef1}\\
& = & \sum_{n = 0}^{\infty}\frac{1}{n!}\langle :\cdot^{\otimes n}:, \xi^{\otimes n}\rangle,
\end{eqnarray}
where $\langle \cdot$ , $\xi \rangle^{\diamond n} : =  \langle \cdot$ , $\xi \rangle \diamond
\langle \cdot$ , $\xi \rangle \diamond \cdots \diamond \langle \cdot$ , $\xi \rangle$
($n$ times). The point--wise formula for $\varphi_{\xi}$ is:
\begin{eqnarray}
\varphi_{\xi}(x) & = & e^{\langle x, \xi \rangle - \frac{1}{2}\langle
\xi, \xi \rangle}, \label{expdef2}
\end{eqnarray}
for almost all $x \in {\mathcal E}'$. \\
It is easy to see that $\varphi_{\xi}$ belongs to
$L^p({\mathcal E}'$, $\mu)$, for all $1 \leq p < \infty$,
and all $\xi$ in $E_c$. The family of exponential functions is
closed with respect to the Wick product and second quantization
operators. Moreover, the vector space spanned by the exponential
functions is closed with respect to each $T$--Wick product, for
every bounded self--adjoint and non--negative operator
$T$, and dense in every space
$L^p({\mathcal E}'$, $\mu)$, for $1 \leq p < \infty$.
We have the lemma:
\begin{lemma}
For all $\xi$ and $\eta$ in $E_c$, and $T$ a bounded self--adjoint operator on $E_c$, such that
$T \geq 0$, we have:
\begin{itemize}
\item
\begin{eqnarray}
\Gamma(T)\varphi_{\xi} & = & \varphi_{T\xi}.
\end{eqnarray}
\item
\begin{eqnarray}
\varphi_{\xi} \diamond_T \varphi_{\eta} & = & e^{\langle T\xi, \eta \rangle}\varphi_{\xi + \eta}.
\end{eqnarray}
\item
\begin{eqnarray}
\varphi_{\xi} \diamond \varphi_{\eta} & = & \varphi_{\xi + \eta}.
\end{eqnarray}
\end{itemize}
\end{lemma}
\begin{proof}
Since for all $n \geq 0$, $\langle :\cdot^{\otimes n}:, \xi^{\otimes n} \rangle \in {\mathcal G}_n$, we have:
\begin{eqnarray*}
\Gamma(T)\varphi_{\xi} & = & \sum_{n = 0}^{\infty}\frac{1}{n!}\langle :\cdot^{\otimes n}:, T^{\otimes n}\xi^{\otimes n} \rangle\\
& = & \sum_{n = 0}^{\infty}\frac{1}{n!}\langle :\cdot^{\otimes n}:, (T\xi)^{\otimes n} \rangle\\
& = & \varphi_{T\xi}.
\end{eqnarray*}
For all $x \in {\mathcal E}'$, we have:
\begin{eqnarray*}
\left(\varphi_{\xi} \diamond_T \varphi_{\eta}\right)(x) & = &
\Gamma\left(\frac{1}{\sqrt{T}}\right)
\left[\Gamma(\sqrt{T})\varphi_{\xi} \cdot \Gamma(\sqrt{T})\varphi_{\eta}\right](x)\\
& = & \Gamma\left(\frac{1}{\sqrt{T}}\right)
\left[\varphi_{\sqrt{T}\xi} \cdot \varphi_{\sqrt{T}\eta}\right](x)\\
& = & \Gamma\left(\frac{1}{\sqrt{T}}\right)
\left[e^{\langle x, \sqrt{T}\xi \rangle - (1/2)\langle\sqrt{T}\xi, \sqrt{T}\xi\rangle}
\cdot e^{\langle x, \sqrt{T}\eta \rangle - (1/2)\langle\sqrt{T}\eta, \sqrt{T}\eta\rangle}\right]\\
& = & \Gamma\left(\frac{1}{\sqrt{T}}\right)
\left[e^{\langle x, \sqrt{T}(\xi + \eta) \rangle - (1/2)(\langle T\xi, \xi\rangle + \langle T\eta, \eta\rangle)}\right].
\end{eqnarray*}
In the exponential we subtract and add $\langle \sqrt{T}\xi$, $\sqrt{T}\eta\rangle$ to complete the square, and obtain:
\begin{eqnarray*}
\left(\varphi_{\xi} \diamond_T \varphi_{\eta}\right)(x) & = & e^{\langle T\xi, \eta\rangle}
\Gamma\left(\frac{1}{\sqrt{T}}\right)
\left[e^{\langle x, \sqrt{T}(\xi + \eta) \rangle - 1/2\langle\sqrt{T}(\xi + \eta), \sqrt{T}(\xi + \eta)\rangle}\right]\\
& = & e^{\langle T\xi, \eta\rangle}
\left[\Gamma\left(\frac{1}{\sqrt{T}}\right)\varphi_{\sqrt{T}(\xi + \eta)}\right]\\
& = & e^{\langle T\xi, \eta\rangle}
\varphi_{(1/\sqrt{T})\sqrt{T}(\xi + \eta)}(x)\\
& = & e^{\langle T\xi, \eta\rangle}\varphi_{\xi + \eta}(x).
\end{eqnarray*}
For $T := 0$, we obtain:
\begin{eqnarray*}
\varphi_{\xi} \diamond \varphi_{\eta} & = & \varphi_{\xi + \eta}.
\end{eqnarray*}
\end{proof}\\
We also have the following functorial property:
\begin{lemma}\label{functorial_lemma}
Let $B$ and $T$ be two commuting, bounded, non--negative, and self--adjoint operators on $E_c$, such that
$B$ is invertible.
Then there exists a vector space $V$, that is dense in all the spaces $L^p({\mathcal E}'$, $\mu)$,
with $1 \leq p < \infty$, such that for any two random variables $\varphi$ and $\psi$ in $V$, we have:
\begin{eqnarray}
\Gamma(B)\left(\varphi \diamond_T \psi\right) & = & \Gamma(B)\varphi \diamond_{TB^{-2}} \Gamma(B)\psi.
\label{functorial}
\end{eqnarray}
\end{lemma}
\begin{proof}
We can take $V$ to be the vector space spanned by the exponential functions. For any $\varphi$ and $\psi$
in $V$, we have:
\begin{eqnarray*}
\Gamma(B)\left(\varphi \diamond_T \psi\right) & = & \Gamma(B)\Gamma\left({\frac{1}{\sqrt{T}}}\right)
\left[\Gamma\left(\sqrt{T}\right)\varphi \cdot \Gamma\left(\sqrt{T}\right)\psi\right]\\
& = & \Gamma\left(\frac{B}{\sqrt{T}}\right)\left[\Gamma\left(\frac{\sqrt{T}}{B}\right)\Gamma(B)\varphi \cdot
\Gamma\left(\frac{\sqrt{T}}{B}\right)\Gamma(B)\psi\right]\\
& = & \Gamma(B)\varphi \diamond_{TB^{-2}} \Gamma(B)\psi.
\end{eqnarray*}
\end{proof}


\section{An integral representation of the generalized Gaussian Wick product}

In this section we prove
an integral representation of the Gaussian $T$--Wick product for every $0 \leq T \leq 2I$, where
$I$ denotes the identity operator.\\
Let us denote by $Exp$ the complex vector space generated by
the exponential functions with subscripts from ${\mathcal E}_c$, that means,
\begin{eqnarray}
Exp & = & \{c_1\varphi_{\xi_1} + \cdots + c_n\varphi_{\xi_n}
\mid n \geq 1, c_i \in {\mathbb C}, \xi_i \in {\mathcal E}_c, 1 \leq i \leq n\}.
\end{eqnarray}
We have the following lemma:
\begin{lemma}\label{lema_de_conectare}
Let $T$ be a self--adjoint and diagonal operator on $E_c$,
commuting with the operator $A$ used in the construction of the Gel'fand triple
${\mathcal E} \subset E \subset {\mathcal E}'$, such that:
\begin{eqnarray}
0 \leq & T & \leq 2I,
\end{eqnarray}
where $I$ denotes the identity operator of $E_c$.\\
Let $C$ and $D$ be two bounded, self--adjoint and diagonal operators on $E_c$, of operatorial norm
at most $1$, commuting with $A$, such that:
\begin{eqnarray}
\left(I - C^2\right)\left(I - D^2\right) & \geq & (T - I)^2C^2D^2.
\end{eqnarray}
That means, if $\{e_n\}_{n \geq 0}$ is the orthonormal basis of $E$ made up of eigenfunctions of the
operator $A$, then there are three sequences
of real numbers, $\{t_n\}_{n \geq 0}$, $\{\alpha_n\}_{n \geq 0}$, and $\{\beta_n\}_{n \geq 0}$, such that:
\begin{itemize}

\item For all $f \in E$, we have:

\begin{eqnarray}
Tf & = & \sum_{n = 0}^{\infty}t_n\langle f, e_n \rangle e_n,
\end{eqnarray}
\begin{eqnarray}
Cf & = & \sum_{n = 0}^{\infty}\alpha_n\langle f, e_n \rangle e_n,
\end{eqnarray}
and
\begin{eqnarray}
Df & = & \sum_{n = 0}^{\infty}\beta_n\langle f, e_n \rangle e_n.
\end{eqnarray}

\item For all $n \geq 0$, $0 \leq t_n \leq 2$, $-1 \leq \alpha_n \leq 1$, and $-1 \leq \beta_n \leq 1$.

\item For all $n \geq 0$, we have:
\begin{eqnarray}
\left(1 - \alpha_n^2\right)\left(1 - \beta_n^2\right) & \geq & (t_n - 1)^2\alpha_n^2\beta_n^2
\label{alpha_beta_t_inequality}.
\end{eqnarray}

\end{itemize}
Then, for all $\varphi$ and $\psi$ in $Exp$, we have:
\begin{eqnarray}
& \ & (\Gamma(C)\varphi \diamond_T \Gamma(D)\psi)(x) \nonumber\\
& = & \int_{{\mathcal E}'}\int_{{\mathcal E}'}\varphi\left(C^*x + P^*y + Q^*z\right)
\psi\left(D^*x + R^*y + S^*z\right)d\mu(z)d\mu(y) \label{integral_representation}\\
& = & \lim_{n \to \infty}\int_{E_n}\int_{E_n}
E\left[\varphi \mid {\mathcal F}_n\right]\left(Cx + Py + Qz\right) \nonumber\\
& \ & \ \ \ \ \ \ \ \ \ \ \ \ \ \ \ \ \ \times
E\left[\psi \mid {\mathcal F}_n\right]\left(Dx + Ry + Sz\right)d\mu_n(z)d\mu_n(y),
\label{finite_integral_representation}
\end{eqnarray}
where $P$, $Q$, $R$, and $S$ are any bounded, self--adjoint, and diagonal operators
commuting with $A$, such that:
\begin{eqnarray}
P^2 + Q^2 & = & I - C^2, \label{Pythagoras1}
\end{eqnarray}
\begin{eqnarray}
R^2 + S^2 & = & I - D^2, \label{Pythagoras2}
\end{eqnarray}
and
\begin{eqnarray}
PR + QS & = & (T - I)CD, \label{produs_scalar}
\end{eqnarray}
and for all $n \in {\mathbb N}$, ${\mathcal F}_n$ denotes the smallest sigma algebra with respect
to which $\langle \cdot$, $e_0 \rangle$, $\langle \cdot$, $e_1 \rangle$, $\cdots$, $\langle \cdot$,
$e_{n - 1} \rangle$ are measurable,
\begin{eqnarray}
E_n & = & {\mathbb R}e_0 \oplus {\mathbb R}e_1 \oplus \cdots \oplus {\mathbb R}e_{n - 1},
\end{eqnarray}
and $\mu_n$ is the standard Gaussian probability measure on $E_n \equiv {\mathbb R}^n$.
The limit from formula (\ref{finite_integral_representation}) is in the $L^p$--sense, for
all $1 \leq p < \infty$.\\
The operators $C^*$, $D^*$, $P^*$, $Q^*$, $R^*$, and $S^*$ map ${\mathcal E}'$ into
${\mathcal E}'$ and represent the dual operators of $C$, $D$, $P$, $Q$, $R$, and $S$ viewed as operators
from ${\mathcal E}$ to ${\mathcal E}$.
\end{lemma}

\begin{proof}
We prove first the existence of such operators $P$, $Q$, $R$ and $S$.\\
For every $n \geq 0$, since $1 - \alpha_n^2 \geq 0$, $1 - \beta_n^2 \geq 0$, and
$(1 - \alpha_n^2)(1 - \beta_n^2) \geq (t_n - 1)^2\alpha_n^2\beta_n^2$, there exist
two vectors $(p_n$, $q_n)$ and $(r_n$, $s_n)$ in ${\mathbb R}^2$, such that:
\begin{eqnarray}
\parallel (p_n, q_n) \parallel_2 & = & \sqrt{1 - \alpha_n^2},
\end{eqnarray}
\begin{eqnarray}
\parallel (r_n, s_n) \parallel_2 & = & \sqrt{1 - \beta_n^2},
\end{eqnarray}
and
\begin{eqnarray}
(p_n, q_n) \cdot (r_n, s_n) & = & (t_n - 1)\alpha_n\beta_n.
\end{eqnarray}
Geometrically, it means that the vectors $(p_n$, $q_n)$ and $(r_n$, $s_n)$ have length
$\sqrt{1 - \alpha_n^2}$ and $\sqrt{1 - \beta_n^2}$, respectively, and the angle between
them have a radian measure of $\arccos((t_n - 1)\alpha_n\beta_n/\sqrt{(1 - \alpha_n^2)(1 - \beta_n^2)})$,
if $|\alpha_n| < 1$ and $|\beta_n| < 1$. If $\alpha_n = \pm 1$ or $\beta_n = \pm 1$, then
(\ref{alpha_beta_t_inequality}) implies $(t_n - 1)\alpha_n\beta_n = 0$ and so we can take
$(p_n$, $q_n) = (0$, $0)$ or $(r_n$, $s_n) = (0$, $0)$.
Once a pair of vectors $(p_n$, $q_n)$ and $(r_n$, $s_n)$ is found, any rotation by an angle
$\theta_n$ of these vectors will produce another pair with the same properties.\\
We choose for each $n \geq 0$, a pair of vectors $(p_n$, $q_n)$ and $(r_n$, $s_n)$ with the above
properties, and define the operators $P$, $Q$, $R$, and $S$ in the following way. For every
$f$ in $E$, let:
\begin{eqnarray}
Pf & := & \sum_{n = 0}^{\infty}p_n\langle f, e_n \rangle e_n,
\end{eqnarray}
\begin{eqnarray}
Qf & := & \sum_{n = 0}^{\infty}q_n\langle f, e_n \rangle e_n,
\end{eqnarray}
\begin{eqnarray}
Rf & := & \sum_{n = 0}^{\infty}r_n\langle f, e_n \rangle e_n,
\end{eqnarray}
and
\begin{eqnarray}
Sf & := & \sum_{n = 0}^{\infty}s_n\langle f, e_n \rangle e_n.
\end{eqnarray}
Of course, $P$, $Q$, $R$, and $S$ can be extended as linear operators from $E_c$ to $E_c$.
Since for all $n \geq 0$, we have:
\begin{eqnarray}
p_n^2 + q_n^2 & = & 1 - \alpha_n^2\\
& \leq & 1
\end{eqnarray}
and
\begin{eqnarray}
r_n^2 + s_n^2 & = & 1 - \beta_n^2\\
& \leq & 1,
\end{eqnarray}
we conclude that the operators $P$, $Q$, $R$, and $S$
are self--adjoint, diagonal, and bounded operators, of operatorial  norm less than
or equal to $1$, on $E$. Moreover, being diagonalized in the same
basis as $A$, they commute with A.\\
Let us observe that $C$, $D$, $P$, $Q$, $R$, and $S$ map ${\mathcal E}$ into ${\mathcal E}$,
and are continuous from ${\mathcal E}$ to ${\mathcal E}$.
Indeed, for all $\xi \in {\mathcal E}$ and $k \geq 0$,
we have, for example:
\begin{eqnarray*}
|C\xi|_k & = & |A^kC\xi|_0\\
& = & |CA^k\xi|_0\\
& \leq & |A^k\xi|_0\\
& = & |\xi|_k.
\end{eqnarray*}
It is clear that relations (\ref{Pythagoras1}), (\ref{Pythagoras2}), and (\ref{produs_scalar})
are satisfied.\\
\par Since both sides of (\ref{integral_representation}) are bilinear in $\varphi$ and $\psi$, to
prove it for every linear combination of exponential functions, it is enough to check it for $\varphi$
and $\psi$ exponential functions.\\
So let $\varphi = \varphi_{\xi}$ and $\psi = \varphi_{\eta}$, for some $\xi$ and $\eta$ in ${\mathcal E}_c$.
We have:
\begin{eqnarray*}
& \ & \int_{{\mathcal E}'}\int_{{\mathcal E}'}\varphi_{\xi}\left(C^*x + P^*y + Q^*z\right)
\varphi_{\eta}\left(D^*x + R^*y + S^*z\right)d\mu(z)d\mu(y)\\
& = & \int_{{\mathcal E}'}\int_{{\mathcal E}'}
\exp\left(\langle C^*x + P^*y + Q^*z, \xi \rangle -
\frac{1}{2}\langle \xi, \xi \rangle\right)\\
& \ & \ \ \ \ \ \ \ \times \exp\left(\langle D^*x + R^*y + S^*z, \eta \rangle -
\frac{1}{2}\langle \eta, \eta \rangle\right)d\mu(z)d\mu(y)\\
& = & \exp\left(\langle x, C\xi + D\eta \rangle - \frac{1}{2}\langle \xi, \xi \rangle -
\frac{1}{2}\langle \eta, \eta \rangle\right)\\
& \ & \int_{{\mathcal E}'}\exp\left(\langle y, P\xi + R\eta \rangle\right)d\mu(y)
\int_{{\mathcal E}'}\exp\left(\langle z, Q\xi + S\eta \rangle\right)d\mu(z)\\
& = & \exp\left(\langle x, C\xi + D\eta \rangle - \frac{1}{2}\langle \xi, \xi \rangle -
\frac{1}{2}\langle \eta, \eta \rangle\right)\\
& \ & \times \int_{{\mathcal E}'}\varphi_{P\xi + R\eta}(y)\exp\left(\frac{1}{2}\langle P\xi + R\eta, P\xi + R\eta \rangle\right)d\mu(y)\\
& \ & \ \ \ \ \times \int_{{\mathcal E}'}\varphi_{Q\xi + S\eta}(z)\exp\left(\frac{1}{2}\langle Q\xi + S\eta, Q\xi + S\eta \rangle\right)d\mu(z).
\end{eqnarray*}
Taking now the constant factors $\exp((1/2)\langle P\xi + R\eta, P\xi + R\eta \rangle)$ and
$\exp((1/2)\langle Q\xi + S\eta, Q\xi + D\eta \rangle)$ outside from their integrals, we obtain:
\begin{eqnarray*}
& \ & \int_{{\mathcal E}'}\int_{{\mathcal E}'}\varphi_{\xi}\left(C^*x + P^*y + Q^*z\right)
\varphi_{\eta}\left(D^*x + R^*y + S^*z\right)d\mu(z)d\mu(y)\\
& = & \exp\left(\langle x, C\xi + D\eta \rangle - \frac{1}{2}\langle (I - P^2 - Q^2)\xi, \xi \rangle\right)\\
& \ & \times \exp\left(\langle (PR + QS)\xi, \eta\rangle -
\frac{1}{2}\langle(I - R^2 - S^2)\eta, \eta \rangle\right)\\
& \ & \times E\left[\varphi_{P\xi + R\eta}\right]E\left[\varphi_{Q\xi + S\eta}\right].
\end{eqnarray*}
Since the expectation of every exponential function is $1$, we have:
\begin{eqnarray*}
& \ & \int_{{\mathcal E}'}\int_{{\mathcal E}'}\varphi_{\xi}\left(C^*x + P^*y + Q^*z\right)
\varphi_{\eta}\left(D^*x + R^*y + S^*z\right)d\mu(z)d\mu(y)\\
& = & \exp\left(\langle x, C\xi + D\eta \rangle - \frac{1}{2}\langle C^2\xi, \xi\rangle + \langle(T - I)CD\xi, \eta\rangle
 - \frac{1}{2}\langle D^2\eta, \eta\rangle \right)\\
& = & \exp(\langle TC\xi, D\eta \rangle)\\
& \ & \times \exp\left(\langle x, C\xi + D\eta \rangle - \frac{1}{2}\langle C\xi, C\xi\rangle - \langle C\xi, D\eta\rangle
 - \frac{1}{2}\langle D\eta, D\eta\rangle \right)\\
& = & \exp(\langle TC\xi, D\eta \rangle)\varphi_{C\xi + D\eta}\\
& = & \varphi_{C\xi} \diamond_T \varphi_{D\eta}\\
& = & \Gamma(C)\varphi_\xi \diamond_T  \Gamma(D)\varphi_\eta.
\end{eqnarray*}
It is not hard to see that, for all $n \geq 1$ and all $u \in E$, we have:
\begin{eqnarray}
E\left[\varphi_u \mid {\mathcal F}_n\right] & = & \varphi_{u_n},
\end{eqnarray}
where:
\begin{eqnarray}
u_n & := & \sum_{k = 0}^{n - 1}\langle u, e_k \rangle e_k.
\end{eqnarray}
Since
\begin{eqnarray}
\lim_{n \to \infty}\varphi_{u_n} = \varphi_u,
\end{eqnarray}
where this limit is computed in the $L^p$--norm, for all $1 \leq p < \infty$,
formula (\ref{finite_integral_representation}) follows easily.
\end{proof}

\begin{corollary}\label{corolar_de_connectare}
Let $T$ be a self--adjoint and diagonal operator on $E_c$ commuting with the operator $A$, such that
$0 \leq T \leq 2I$.\\
Let $C$ and $D$ be two bounded, self--adjoint and diagonal operators on $E_c$, of operatorial norm
at most $1$, commuting with $A$, such that:
\begin{eqnarray}
\left(I - C^2\right)\left(I - D^2\right) & = & (T - I)^2C^2D^2.
\end{eqnarray}
Then, for all $\varphi$ and $\psi$ in $Exp$, we have:
\begin{eqnarray}
& \ & (\Gamma(C)\varphi \diamond_T \Gamma(D)\psi)(x) \nonumber\\
& = & \int_{{\mathcal E}'}\varphi\left(C^*x + P^*y\right)
\psi\left(D^*x + R^*y\right)d\mu(y), \label{simpler_integral_representation}
\end{eqnarray}
where
\begin{eqnarray}
P & = & {\rm sgn}\left[(T - I)CD\right]\sqrt{I - C^2} \label{Pytha1}
\end{eqnarray}
and
\begin{eqnarray}
R & = & \sqrt{I - D^2}, \label{Pytha2}
\end{eqnarray}
where ${\rm sgn}$ is the right--continuous extension, to ${\mathbb R}$, of
the function $x \mapsto x/|x|$, defined on ${\mathbb R} \setminus \{0\}$.
Here ${\rm sgn}[(T - I)CD]$ is computed according to the Dunford functional calculus,
which in this case, since $T - I$, $C$, and $D$ are self--adjoint and diagonal operators,
that are diagonalized in the same basis, means to simply apply the measurable function
${\rm sgn}$ to each eigenvalue of the operator $(T - I)CD$.

\end{corollary}

\begin{proof}
In Lemma \ref{lema_de_conectare} take $Q = S := 0$,
$P = {\rm sgn}\left[(T - I)CD\right]\sqrt{I - C^2}$ and
$R = \sqrt{I - D^2}$.
\end{proof}

\begin{observation}
Lemma \ref{lema_de_conectare} remains true if the condition that each of the operators
$C$, $D$, and $T$ commutes with $A$, is replaced by the fact that $C$, $D$, and $T$
commute among themselves and are continuous (bounded) linear operators from ${\mathcal E}$ to
${\mathcal E}$. In particular, in the finite dimensional case (that means, if the dimension of $E$
is a finite number $d$, in which case ${\mathcal E} = E = {\mathcal E}' \equiv {\mathbb R}^d$),
since the condition of being bounded is automatically satisfied, the only condition required
is that $C$, $D$, and $T$ commute among themselves.
\end{observation}

The above observation is important for people who are not interested in the infinite dimensional
case, and are content with the finite dimensional one. In that case, the technicality of commuting
with the operator $A$ is removed. In fact, in that case there is no need to speak of such an operator
$A$, since a Gaussian probability measure on a finite dimensional space can be defined without it.

\section{H\"older inequalities for norms of generalized Wick products}

We present now the main result of this paper. To prove our result we need
the following theorem of Lieb (see \cite{l1990}
or \cite{ll2001} (page 100)).

\begin{theorem} \label{Lieb}
Fix $k > 1$, integers $n_1$, $\dots$, $n_k$ and numbers $p_1$,
$\dots$, $p_k \geq 1$. Let $M \geq 1$ and let $B_i$ (for $i = 1$,
$\dots$, $k$) be a linear mapping from ${\mathbb R}^M$ to ${\mathbb
R}^{n_i}$. Let $Z : {\mathbb R}^M \to {\mathbb R}^+$ be some fixed
Gaussian function,
\begin{eqnarray*}
Z(x) & = & \exp\left[-\langle x, Jx \rangle\right]
\end{eqnarray*}
with $J$ a real, positive semi--definite $M \times M$ matrix
(possible zero).
\par For functions $f_i$ in $L^{p_i}({\mathbb R}^{n_i})$ consider the
integral $I_Z$ and the number $C_Z$:
\begin{eqnarray}
I_Z(f_1, \dots, f_k) & = & \int_{{\mathbb R}^M}Z(x)\prod_{i =
1}^kf_i(B_ix)dx, \label{I_Z}
\end{eqnarray}
\begin{eqnarray}
C_Z & := & \sup\{I_Z(f_1, \dots, f_k) \mid \||f_i\||_{p_i} = 1 \ for
\ i = 1, \dots, k\}, \label{C_Z}
\end{eqnarray}
where $|\| \cdot \||_{p_i}$ denotes $L^{p_i}({\mathbb R}^{n_i})$ norm
with respect to the
Lebesgue measure.
Then $C_Z$ is determined by restricting the $f$'s to be Gaussian
functions, i.e.,
\begin{eqnarray*}
C_Z & = & \sup\{I_Z(f_1, \dots, f_k) \mid \||f_i\||_{p_i} = 1 \ and
\ f_i(x) =
c_i\exp[-\langle x, J_ix \rangle]\\
& \ & with \ c_i > 0, \ and \ J_i \ a \ real, \ symmetric, \
positive-definite\\
& \ & \ n_i \times n_i \ matrix\}.
\end{eqnarray*}
\end{theorem}
The proof of the next corollary, of the above theorem, can be found
in \cite{dpls11}.
\begin{corollary} \label{corlieb}
Let $p$, $q$, $r \geq 1$. Let $B_1$ and $B_2$ be linear maps from
${\mathbb R}^2$ to ${\mathbb R}^2$, and $J$ a real,
positive--semidefinite $2 \times 2$ matrix (possible zero). For $f$
in $L^p({\mathbb R}^2)$ and $g$ in $L^q({\mathbb R}^2)$, we consider
the product:
\begin{eqnarray}
\left(f \star_{B_1, B_2, J} g\right)(x) & = & \int_{{\mathbb
R}}f(B_1(x, y))g(B_2(x, y))e^{-\langle (x, y), J(x, y)\rangle}dy.
\end{eqnarray}
We define:
\begin{eqnarray}
C & := & \sup\{|\| f \star_{B_1, B_2, J} g \||_r \mid |\|f \||_p =
|\| g \||_q = 1\}.
\end{eqnarray}
Then $C$ is determined by restricting $f$ and $g$ to be Gaussian
functions.
\end{corollary}


We will also need the following computational lemma:

\begin{lemma}\label{computational_lemma}
Let $\{a_i\}_{1 \leq i \leq n}$ and $\{b_i\}_{1 \leq i \leq n}$
be two finite sequences of real numbers. Then, we have:
\begin{eqnarray}
& \ & \frac{1}{\sqrt{2\pi}}\int_{{\mathbb R}}\exp\left[-\frac{1}{2}\sum_{i = 1}^n(a_ix + b_i)^2\right]dx \nonumber\\
& = & \frac{1}{\sqrt{\sum_{i = 1}^na_i^2}}\exp\left[-\frac{1}{2}\frac{\sum_{i = 1}^na_i^2\sum_{i = 1}^nb_i^2 -
\left(\sum_{i = 1}^na_ib_i\right)^2}{\sum_{i = 1}^na_i^2}\right] \label{Cauchy_Buniakovski_Schwarz}\\
& = & \frac{1}{\sqrt{\sum_{i = 1}^na_i^2}}\exp\left[-\frac{1}{2}\frac{\sum_{1 \leq i < j \leq n}
\left(a_ib_j - a_jb_i\right)^2}{\sum_{i = 1}^na_i^2}\right]. \label{Lagrange}
\end{eqnarray}

\begin{proof}
We have:
\begin{eqnarray*}
& \ & \frac{1}{\sqrt{2\pi}}\int_{{\mathbb R}}\exp\left[-\frac{1}{2}\sum_{i = 1}^n(a_ix + b_i)^2\right]dx\\
& = & \exp\left[-\frac{1}{2}\sum_{i = 1}^nb_i^2\right]\\
& \ & \times \frac{1}{\sqrt{2\pi}}\int_{{\mathbb R}}
\exp\left[-\frac{1}{2}\left(\sum_{i = 1}^na_i^2\right)x^2 + \left(\sum_{i = 1}^na_ib_i\right)x\right]dx
\end{eqnarray*}
Making first the substitution $x' = x\sqrt{\sum_{i = 1}^na_i^2}$ and then completing the square in the
exponential, we obtain:\\
\begin{eqnarray*}
& \ & \frac{1}{\sqrt{2\pi}}\int_{{\mathbb R}}\exp\left[-\frac{1}{2}\sum_{i = 1}^n(a_ix + b_i)^2\right]dx\\
& = & \frac{\exp\left[-(1/2)\sum_{i = 1}^nb_i^2\right]}{\sqrt{\sum_{i = 1}^na_i^2}}
\cdot \frac{1}{\sqrt{2\pi}}\int_{{\mathbb R}}
\exp\left[-\frac{1}{2}x'^2 + \frac{\sum_{i = 1}^na_ib_i}{\sqrt{\sum_{i = 1}^na_i^2}}x'\right]dx'\\
& = & \frac{\exp\left[-(1/2)\sum_{i = 1}^nb_i^2\right]}{\sqrt{\sum_{i = 1}^na_i^2}}
\exp\left[+\frac{1}{2}\frac{\left(\sum_{i = 1}^na_ib_i\right)^2}{\sum_{i = 1}^na_i^2}\right]\\
& \ & \times \frac{1}{\sqrt{2\pi}}\int_{{\mathbb R}}
\exp\left[-\frac{1}{2}\left(x' - \frac{\sum_{i = 1}^na_ib_i}{\sqrt{\sum_{i = 1}^na_i^2}}\right)^2\right]dx'\\
& = &  \frac{1}{\sqrt{\sum_{i = 1}^na_i^2}}\exp\left[-\frac{1}{2}\frac{\sum_{i = 1}^na_i^2\sum_{i = 1}^nb_i^2 -
\left(\sum_{i = 1}^na_ib_i\right)^2}{\sum_{i = 1}^na_i^2}\right].
\end{eqnarray*}
Formula (\ref{Lagrange}) follows now from Lagrange identity.
\end{proof}

\end{lemma}

\begin{theorem} \label{holder_young_theorem} {\bf (H\"older--Young
inequality for generalized Gaussian Wick products)} Let
${\mathcal E} \subset E \subset {\mathcal E}'$ be a Gel'fand triple
given by a self--adjoint diagonal operator $A$ on $E$, with increasing, greater
than $1$ eigenvalues, whose inverse is a Hilbert--Schmidt operator. Let
$\mu$ be the Gaussian probability measure on ${\mathcal E}'$ whose existence
is guaranteed by Minlos Theorem.
Let $T$ be a self--adjoint, diagonal operator on $E$, commuting with $A$, such that:
\begin{eqnarray}
0 & \leq~ T ~\leq & 2I,
\end{eqnarray}
where $I$ denotes the identity operator of $E$.
Let $C$ and $D$ be two self--adjoint, diagonal, and bounded operators
on $E$, of operatorial norm less than or equal to $1$, commuting with the
operator $A$, such that:
\begin{eqnarray}
\left(I - C^2\right)(I - D^2) & \geq & (T - I)^2C^2D^2. \label{C_D_t_condition}
\end{eqnarray}
Let $p$, $q$, $r > 1$ such that:\\
\begin{eqnarray}
(r - 1)I & \leq & \frac{(p - 1)(q - 1)I - C^2D^2T^2}
{(q - 1)C^2 + (p - 1)D^2 + 2C^2D^2T}
\label{Holder_condition}
\end{eqnarray}
or equivalently:
\begin{eqnarray}
\frac{1}{(r - 1)I + T} & \geq & \frac{1}{\frac{p - 1}{C^2} + T} +
\frac{1}{\frac{q - 1}{D^2} + T}, \label{equivalent_Holder_condition}
\end{eqnarray}
with the convention that if $x \in Ker(C) := C^{-1}(0)$, then the first fraction operator
from the right--hand side of (\ref{equivalent_Holder_condition}) evaluated at $x$ is zero,
and a similar convention
for $x \in Ker(D)$.
Then for all $\varphi$ in $L^p({\mathcal E}'$, $\mu)$ and
$\psi$ in $L^q({\mathcal E}'$, $\mu)$,
$\Gamma(C)\varphi \diamond_T \Gamma(D)\psi$ belongs
to $L^r({\mathcal E}'$, $\mu)$, and the following
inequality holds:
\begin{eqnarray}
\left\|\Gamma(C)\varphi \diamond_T
\Gamma(D)\psi\right\|_r & \leq & \|\varphi\|_p \cdot
\|\psi\|_q. \label{fullholder}
\end{eqnarray}
On the other hand, if:
\begin{eqnarray*}
(r - 1)I & \not\leq & \frac{(p - 1)(q - 1)I - C^2D^2T^2}
{(q - 1)C^2 + (p - 1)D^2 + 2C^2D^2T},
\end{eqnarray*}
then the bilinear operator $(\varphi$, $\psi) \mapsto \Gamma(C)\varphi \diamond_T \Gamma(D)\psi$
is not bounded from  $L^p({\mathcal E}'$, $\mu) \times L^q({\mathcal E}'$, $\mu)$ to  $L^r({\mathcal E}'$, $\mu)$.
\end{theorem}


We present now the proof of Theorem \ref{holder_young_theorem}. The proof will be subdivided into many steps.
The main idea of each step will be written in Italic.\\

\begin{proof}\\
($\Rightarrow$) Let us assume first that:
\begin{eqnarray*}
(r - 1)I & \leq & \frac{(p - 1)(q - 1)I - C^2D^2T^2}
{(q - 1)C^2 + (p - 1)D^2 + 2C^2D^2T}.
\end{eqnarray*}
{\bf Step 0.} \ {\em We may assume that $T > 0$, $C^2 > 0$, and $D^2 > 0$.}\\
\ \\
Indeed, since $T \geq 0$ can be obtained from the case $T' > 0$, by
gently passing to the limit on $Ker(T) = T^{-1}(0)$ as $P_{Ker(T)}T' \to 0^+$, where
$P_{Ker(T)}$ denotes the orthogonal projection of $E$ on $Ker(T)$,
we may assume that $T > 0$. Similarly, we may assume that $C^2 > 0$ and $D^2 > 0$.\\
\ \\
{\bf Step 1.} \ {\em In the finite dimensional case, we can reduce the problem to
 the one--dimensional case via Minkowski integral inequality.}\\
\ \\
Let $d$ be a fixed finite dimension and let $\mu_d$ be the standard Gaussian probability measure on ${\mathbb R}^d$.
Let $T$ be in $(0$, $2I]$, and $C$ and $D$ be two non--zero commuting linear self--adjoint operators from ${\mathbb R}^d$ to
${\mathbb R}^d$, and $p$, $q$, $r \geq 1$, such that
relations (\ref{C_D_t_condition}) and (\ref{Holder_condition}) hold. Let $\{e_i\}_{1 \leq i \leq d}$ be an orthogonal
basis of ${\mathbb R}^d$ that diagonalizes all three
operators $T$, $C$ and $D$. If $x = x_1e_1 + x_2e_2 + \cdots + x_de_d \in {\mathbb R}^d$, then we write $x = (x_1$, $x_2$,
$\dots$, $x_d)$.\\
Let $\varphi$ and $\psi$ be two random variables
in $Exp$ (that means, $\varphi$ and $\psi$ are linear combinations of exponential functions). From Lemma \ref{lema_de_conectare},
we know that:
\begin{eqnarray*}
& \ & (\Gamma(C)\varphi \diamond_T \Gamma(D)\psi)(x) \nonumber\\
& = & \int_{{\mathbb R}^d}\int_{{\mathbb R}^d}\varphi\left(C^*x + P^*y + Q^*z\right)
\psi\left(D^*x + R^*y + S^*z\right)d\mu_d(z)d\mu_d(y),
\end{eqnarray*}
where $P$, $Q$, $R$, and $S$ are linear self--adjoint operators on ${\mathbb R}^d$ commuting with $C$ and $D$, such that:
\begin{eqnarray*}
P^2 + Q^2 & = & I - C^2,
\end{eqnarray*}
\begin{eqnarray*}
R^2 + S^2 & = & I - D^2,
\end{eqnarray*}
and
\begin{eqnarray*}
PR + QS & = & (T - I)CD.
\end{eqnarray*}
In fact, in this finite--dimensional case the use of the adjoint ``$\ ^*$" \ is not necessary, because $C = C^*$, $D = D^*$, $\dots$, $S = S^*$.\\
Let us assume that for a fixed triplet, $(p$, $q$, $r)$, inequality (\ref{fullholder}) holds for two finite dimensions $d = m$ and $d = n$, and
prove that it continues to hold for
$d = m + n$. Expanding the vectors $x$, $y$, and $z$ along the orthonormal basis $\{e_i\}_{1 \leq i \leq m + n}$ we can write
$x = (x_m$, $x_n)$, $y = (y_m$, $y_n)$, and $z = (z_m$, $z_n)$, where $(x_m$, $0)$ denotes the projection of $x$
onto the vector space spanned by $\{e_i\}_{1 \leq i \leq m}$, and $(0$, $x_n) = x - (x_m$, $0)$. We also denote by $C_m$ and $C_n$,
the restriction of $C$ to the spaces spanned by $\{e_i\}_{1 \leq i \leq m}$ and $\{e_i\}_{m + 1 \leq i \leq m + n}$, respectively.
Using Minkowski integral inequality we have:
\begin{eqnarray*}
& \ & \parallel \Gamma(C)\varphi \diamond_T \Gamma(D)\psi \parallel_r\\
& = & \left\{\int_{{\mathbb R}^{m + n}}\left|\int_{{\mathbb R}^{m + n}}\int_{{\mathbb R}^{m + n}}
\varphi(Cx + Py + Qz)\psi(Dx + Ry + Sz)\right.\right.\\
& \ & \left.\left.d\mu_{m + n}(z)d\mu_{m + n}(y)\right|^rd\mu_{m + n}(x)\right\}^{1/r}\\
& \leq & \left\{\int_{{\mathbb R}^m}\left\{\int_{{\mathbb R}^n}
\left|\int_{{\mathbb R}^m}\int_{{\mathbb R}^m}\int_{{\mathbb R}^n}\int_{{\mathbb R}^n}\right.\right.\right.\\
& \ & \left|\varphi\left((C_mx_m, C_nx_n)  + (P_my_m, P_ny_n) + (Q_mz_m, Q_nz_n)\right)\right|\\
& \ & \left|\psi\left((D_mx_m, D_nx_n)  + (R_my_m, R_ny_n) + (S_mz_m, S_nz_n)\right)\right|\\
& \ & \left.\left.\left.d\mu_n(z_n)d\mu_n(y_n)d\mu_m(z_m)d\mu_m(y_m)\right|^rd\mu_n(x_n)\right\}d\mu_m(x_m)\right\}^{1/r}\\
& \leq & \left\{\int_{{\mathbb R}^m}\left\{\int_{{\mathbb R}^m}\int_{{\mathbb R}^m}
\left\{\int_{{\mathbb R}^n}\left|\int_{{\mathbb R}^n}\int_{{\mathbb R}^n}\right.\right.\right.\right.\\
& \ & \left|\varphi\left((C_mx_m, C_nx_n)  + (P_my_m, P_ny_n) + (Q_mz_m, Q_nz_n)\right)\right|\\
& \ & \left|\psi\left((D_mx_m, D_nx_n)  + (R_my_m, R_ny_n) + (S_mz_m, S_nz_n)\right)\right|\\
& \ & \left.\left.\left.\left.d\mu_n(z_n)d\mu_n(y_n)\right|^rd\mu_n(x_n)\right\}^{1/r}d\mu_m(z_m)d\mu_m(y_m)\right\}^rd\mu_m(x_m)\right\}^{1/r}\\
& \leq & \left\{\int_{{\mathbb R}^m}\left\{\int_{{\mathbb R}^m}\int_{{\mathbb R}^m}
\parallel \varphi\left((C_mx_m, \cdot)  + (P_my_m, \cdot) + (Q_mz_m, \cdot)\right)\parallel_p\right.\right.\\
& \ &
\left.\parallel \psi\left((D_mx_m, \cdot)  + (R_my_m, \cdot) + (S_mz_m, \cdot)\right)\parallel_qd\mu_m(z_m)d\mu_m(y_m)\right\}^r\\
& \ & \left.d\mu_m(x_m)\right\}^{1/r}\\
& \leq & \parallel \varphi \parallel_p\parallel \psi \parallel_q.
\end{eqnarray*}
\ \\
\noindent {\bf Step 2.} \ {\em In dimension $d = 1$, we may assume that both inequalities (\ref{C_D_t_condition}) and (\ref{Holder_condition})
are equalities.}\\
\ \\
Indeed, if $d = 1$, then $T = tI$, for some $t \in (0$, $2]$, $C = \alpha I$ and $D = \beta I$, for some $\alpha$ and
$\beta$ in $[-1$, $1] \setminus \{0\}$. Inequality (\ref{C_D_t_condition}) becomes:
\begin{eqnarray}
\left(1 - \alpha^2\right)\left(1 - \beta^2\right) & \geq & \alpha^2\beta^2(t - 1)^2. \label{a_b_t_condition}
\end{eqnarray}
If inequality (\ref{a_b_t_condition}) is strict, then since:
\begin{eqnarray*}
\left[1 - (\pm 1)^2\right]\left(1 - \beta^2\right) & \leq & (\pm 1)^2\beta^2(t - 1)^2,
\end{eqnarray*}
there exists $\alpha_0$ having the same sign as $\alpha$, with
$|\alpha_0| \in (|\alpha|$, $1]$, such that:
\begin{eqnarray}
\left(1 - \alpha_0^2\right)\left(1 - \beta^2\right) & = & \alpha_0^2\beta^2(t - 1)^2.
\end{eqnarray}
Let us suppose that inequality (\ref{fullholder}) holds for the triplet $(\alpha_0$, $\beta$, $t)$ and all
functions $\varphi \in L^p({\mathcal E}'$, $\mu)$ and $\psi \in L^q({\mathcal E}'$, $\mu)$. Then since
$\Gamma((\alpha/\alpha_0)I)$ is a bounded operator from $L^p({\mathcal E}'$, $\mu)$ to $L^p({\mathcal E}'$, $\mu)$
of operatorial norm equal to $1$, we have:
\begin{eqnarray*}
\left\|\Gamma(\alpha I)\varphi \diamond_T
\Gamma(\beta I)\psi\right\|_r & = & \left\|\Gamma(\alpha_0 I)\Gamma\left(\frac{\alpha}{\alpha_0}I\right)\varphi \diamond_T
\Gamma(\beta I)\psi\right\|_r\\
& \leq & \left\|\Gamma\left(\frac{\alpha}{\alpha_0}I\right)\varphi\right\|_p \cdot \|\psi\|_q\\
& \leq & \|\varphi\|_p \cdot
\|\psi\|_q.
\end{eqnarray*}
Also inequality (\ref{Holder_condition}), in dimension $d = 1$, becomes:
\begin{eqnarray}
r - 1 & \leq & \frac{(p - 1)(q - 1) - \alpha^2\beta^2t^2}
{(q - 1)\alpha^2 + (p - 1)\beta^2 + 2\alpha^2\beta^2t}. \label{p_q_r_condition}
\end{eqnarray}
If inequality (\ref{p_q_r_condition}) is strict, then since for $p$ and $q$ fixed,
$\lim_{r \to \infty}(r - 1) = \infty$, there exists $r_0 \in (r$, $\infty)$, such that:
\begin{eqnarray*}
r_0 - 1 & = & \frac{(p - 1)(q - 1) - \alpha^2\beta^2t^2}
{(q - 1)\alpha^2 + (p - 1)\beta^2 + 2\alpha^2\beta^2t}.
\end{eqnarray*}
Suppose inequality $(\ref{fullholder})$ holds for the triplet $(p$, $q$, $r_0)$ and all
$\varphi \in L^p({\mathcal E}'$, $\mu)$ and $\psi \in L^q({\mathcal E}'$, $\mu)$.
Then using Lyapunov inequality, we have:
\begin{eqnarray*}
\left\|\Gamma(\alpha I)\varphi \diamond_T
\Gamma(\beta I)\psi\right\|_r & \leq & \left\|\Gamma(\alpha I)\varphi \diamond_T
\Gamma(\beta I)\psi\right\|_{r_0}\\
& \leq & \|\varphi\|_p \cdot \|\psi\|_q.
\end{eqnarray*}

\noindent {\bf Step 3.} \ {\em Change the problem of proving an inequality about the $p$, $q$, and $r$ norms with
respect to the standard Gaussian probability measure into a problem of proving an inequality about the
$p$, $q$, and $r$ norms with respect to the Lebesgue measure.}\\
\ \\
We assume now that we are in dimension $d = 1$, $C = \alpha I$, $D = \beta I$, $T = tI$,
$-1 \leq \alpha \leq 1$, $\alpha \neq 0$, $-1 \leq \beta \leq 1$, $\beta \neq 0$, $t \in (0$, $2]$,
and we have the following equalities:
\begin{eqnarray}
\left(1 - \alpha^2\right)\left(1 - \beta^2\right) & = & \alpha^2\beta^2(t - 1)^2 \label{a_b_t_equality}
\end{eqnarray}
and
\begin{eqnarray}
r - 1 & = & \frac{(p - 1)(q - 1) - \alpha^2\beta^2t^2}
{(q - 1)\alpha^2 + (p - 1)\beta^2 + 2\alpha^2\beta^2t}. \label{p_q_r_equality}
\end{eqnarray}
Formula (\ref{p_q_r_equality}) is equivalent, via formula (\ref{a_b_t_equality}), to:
\begin{eqnarray}
r & = & pq\frac{1 - \frac{1 - \alpha^2}{p} - \frac{1 - \beta^2}{q}}{(q - 1)\alpha^2 + (p - 1)\beta^2 + 2\alpha^2\beta^2t}.
\label{r_p_q_c_d_equality}
\end{eqnarray}
According to Corollary \ref{corolar_de_connectare}, for every $\varphi$ and $\psi$ that are linear combinations of
exponential functions, we have:
\begin{eqnarray}
& \ & \left[\Gamma(\alpha I)\varphi \diamond_t \Gamma(\beta I)\psi\right](x) \label{a_b_t_conectare}\\
& = & \int_{{\mathbb R}}\varphi\left(\alpha x +
{\rm sgn}[(t - 1)\alpha\beta]\sqrt{1 - \alpha^2}y\right)\psi\left(\beta x + \sqrt{1 - \beta^2}y\right)e^{-y^2/2}dy.
\nonumber
\end{eqnarray}
Let us define the numbers:
\begin{eqnarray}
\gamma & := & {\rm sgn}[\alpha\beta(t - 1)]\sqrt{1 - \alpha^2}
\end{eqnarray}
and
\begin{eqnarray}
\delta & := & \sqrt{1 - \beta^2}.
\end{eqnarray}
Then we have:
\begin{eqnarray}
\alpha^2 + \gamma^2 & = & 1,
\end{eqnarray}
\begin{eqnarray}
\beta^2 + \delta^2 & = & 1,
\end{eqnarray}
and
\begin{eqnarray}
\gamma\delta & = & \alpha\beta(t - 1).
\end{eqnarray}
We observe that, for any $0 < u < \infty$, a measurable function $f(x)$ belongs to $L^u({\mathbb R}$, $\mu)$
if and only if $f(x)e^{-x^2/(2u)}$ belongs to $L^u({\mathbb R}$, $d_Nx)$, where $d_Nx$ denotes the normalized
Lebesgue measure $(1/\sqrt{2\pi})dx$ on ${\mathbb R}$. With this in mind,
we are preparing now to cross the bridge from $L^p({\mathbb R}$, $\mu)$ to $L^p({\mathbb R}$, $d_Nx)$, from
$L^q({\mathbb R}$, $\mu)$ to $L^q({\mathbb R}$, $d_Nx)$, and from $L^r({\mathbb R}$, $\mu)$ to
$L^r({\mathbb R}$, $d_Nx)$. To do this, we
multiply both sides of formula (\ref{a_b_t_conectare}) by $e^{-x^2/(2r)}$ and rewrite the expression inside the integral in
the following way:
\begin{eqnarray*}
& \ & (\Gamma(\alpha I)\varphi \diamond_t \Gamma(\beta I)\psi)(x)e^{-x^2/(2r)}\\
& = & \int_{{\mathbb R}}\varphi(\alpha x + \gamma y)
e^{-(\alpha x + \gamma y)^2/(2p)}\psi(\beta x + \delta y)e^{-(\beta x + \delta y)^2/(2q)}\\
& \ & \ \ \ \ \times e^{(\alpha x + \gamma y)^2/(2p)}e^{(\beta x + \delta y)^2/(2q)}
e^{-y^2/2}e^{-x^2/(2r)}d_Ny.
\end{eqnarray*}
Let us define now the following functions:
\begin{eqnarray}
f(x) & : = & \varphi(x) \cdot e^{-x^2/(2p)}
\end{eqnarray}
and
\begin{eqnarray}
g(x) & : = & \psi(x) \cdot e^{-x^2/(2q)}.
\end{eqnarray}
In this way, our problem of proving that the bilinear operator, which is ``hopefully" defined as:
\begin{eqnarray}
B : L^p({\mathbb R}, \mu) \times L^q({\mathbb R}, \mu) & \to & L^r({\mathbb R}, \mu),
\end{eqnarray}
(we said ``hopefully" since we do not know whether it maps $L^p \times L^q$ into $L^r$),
\begin{eqnarray}
B(\varphi, \psi)(x) & := & \int_{{\mathbb R}}\varphi\left(\alpha x + \gamma y\right)
\psi\left(\beta x + \delta y\right)e^{-y^2/2}d_Ny,
\end{eqnarray}
is a bounded operator of operatorial norm equal to $1$,
becomes the equivalent problem of showing  that the bilinear operator, which is ``hopefully" defined as:
\begin{eqnarray}
\tilde{B} : L^p({\mathbb R}, d_Nx) \times L^q({\mathbb R}, d_Nx) & \to & L^r({\mathbb R}, d_Nx),
\end{eqnarray}
\begin{eqnarray}
\tilde{B}(f, g)(x) \nonumber
& := & \int_{{\mathbb R}}f(\alpha x + \gamma y)
g(\beta x + \delta y) \nonumber\\
& \ & \times e^{(\alpha x + \gamma y)^2/(2p)}e^{(\beta x + \delta y)^2/(2q)}
e^{-y^2/2}e^{-x^2/(2r)}d_Ny,
\end{eqnarray}
is a bounded operator of operatorial norm equal to $1$.\\
Let us define the kernel:
\begin{eqnarray}
J_d(x, y) & := &  e^{(\alpha x + \gamma y)^2/(2p)}e^{(\beta x + \delta y)^2/(2q)}
e^{-y^2/2}e^{-x^2/(2r)}.
\end{eqnarray}
We can write $J_d(x$, $y)$ as an exponential of a quadratic form of $(x$, $y)$, in the following way:
\begin{eqnarray}
J_d(x, y) & = & e^{-(1/2)ax^2 + bxy - (1/2)cy^2},
\end{eqnarray}
where:
\begin{eqnarray}
a & := & \frac{1}{r} - \frac{\alpha^2}{p} - \frac{\beta^2}{q},
\end{eqnarray}
\begin{eqnarray}
b & := & \frac{\alpha\gamma}{p} + \frac{\beta\delta}{q},
\end{eqnarray}
and
\begin{eqnarray}
c & := & 1 - \frac{\gamma^2}{p} - \frac{\delta^2}{q}.
\end{eqnarray}
We make the following observations:
\begin{enumerate}
\item $c > 0$.\\
Indeed, from formula (\ref{r_p_q_c_d_equality}), we can see that:
\begin{eqnarray*}
c & = & \frac{r}{pq}\left[(q - 1)\alpha^2 + (p - 1)\beta^2 + 2\alpha^2\beta^2t\right]\\
& > & 0.
\end{eqnarray*}

\item $b^2 = ac$.\\
Indeed, we have:
\begin{eqnarray*}
ac & = & \left(\frac{1}{r} - \frac{\alpha^2}{p} - \frac{\beta^2}{q}\right)
\left(1 - \frac{1 - \alpha^2}{p} - \frac{1 - \beta^2}{q}\right)\\
& = & \frac{1}{r}\left(1 - \frac{1 - \alpha^2}{p} - \frac{1 - \beta^2}{q}\right)\\
& \ & - \frac{\alpha^2}{p} - \frac{\beta^2}{q} + \frac{\alpha^2(1 - \beta^2) + \beta^2(1 - \alpha^2)}{pq}\\
& \ & + \frac{\alpha^2(1 - \alpha^2)}{p^2} + \frac{\beta^2(1 - \beta^2)}{q^2}.
\end{eqnarray*}
Substituting now $r$ by the right--hand side of formula (\ref{r_p_q_c_d_equality}), we obtain:
\begin{eqnarray*}
ac & = & \frac{\alpha^2(q - 1) + \beta^2(p - 1) + 2\alpha^2\beta^2t}{pq - (1 - \alpha^2)q - (1 - \beta^2)p}
\left(1 - \frac{1 - \alpha^2}{p} - \frac{1 - \beta^2}{q}\right)\\
& \ & - \frac{\alpha^2}{p} - \frac{\beta^2}{q} + \frac{\alpha^2(1 - \beta^2) + \beta^2(1 - \alpha^2)}{pq}\\
& \ & + \frac{\alpha^2(1 - \alpha^2)}{p^2} + \frac{\beta^2(1 - \beta^2)}{q^2}\\
& = & \frac{\alpha^2(q - 1) + \beta^2(p - 1) + 2\alpha^2\beta^2t}{pq}\\
& \ & +  \frac{-\alpha^2q - \beta^2p + \alpha^2(1 - \beta^2) + \beta^2(1 - \alpha^2)}{pq}\\
& \ & + \frac{\alpha^2(1 - \alpha^2)}{p^2} + \frac{\beta^2(1 - \beta^2)}{q^2}\\
& = & \frac{2\alpha\beta\left[\alpha\beta(t - 1)\right]}{pq}\\
& \ & + \frac{\alpha^2(1 - \alpha^2)}{p^2} + \frac{\beta^2(1 - \beta^2)}{q^2}.
\end{eqnarray*}
Since:
\begin{eqnarray*}
\left(1 - \alpha^2\right)\left(1 - \beta^2\right) & = & \alpha^2\beta^2(t - 1)^2,
\end{eqnarray*}
we have:
\begin{eqnarray*}
\alpha\beta(t - 1) & = & {\rm sgn}(\alpha\beta(t - 1))\sqrt{\left(1 - \alpha^2\right)\left(1 - \beta^2\right)}\\
& = & \gamma\delta.
\end{eqnarray*}
Thus we obtain:
\begin{eqnarray*}
ac & = & \frac{2\alpha\beta\gamma\delta}{pq} + \frac{\alpha^2(1 - \alpha^2)}{p^2} + \frac{\beta^2(1 - \beta^2)}{q^2}\\
& = & \frac{2\alpha\beta\gamma\delta}{pq} + \frac{\alpha^2\gamma^2}{p^2} + \frac{\beta^2\delta^2}{q^2}\\
& = & \left(\frac{\alpha\gamma}{p} + \frac{\beta\delta}{q}\right)^2\\
& = & b^2.
\end{eqnarray*}

\item $a \geq 0$.\\
Since we already know that $c > 0$ and:
\begin{eqnarray*}
ac & = & b^2\\
& \geq & 0,
\end{eqnarray*}
we conclude that $a \geq 0$.

\end{enumerate}

The three conditions $a \geq 0$, $b^2 = ac$, and $c > 0$, imply that the quadratic form
$-(1/2)ax^2 + bxy - (1/2)cy^2$
is negative semi--definite. More precisely, if we define:
\begin{eqnarray}
m & := & \sqrt{a}
\end{eqnarray}
and
\begin{eqnarray}
n & := & {\rm sgn}(b)\sqrt{c},
\end{eqnarray}
then we have:
\begin{eqnarray}
-\frac{1}{2}ax^2 + bxy - \frac{1}{2}cy^2 & = &
-\frac{1}{2}(mx - ny)^2.
\end{eqnarray}
For all $u \geq 1$, we will denote the $L^u$--norm with respect to the standard Gaussian measure $\mu$ on ${\mathbb R}$
by $\| \cdot \|_u$, and the $L^u$--norm with respect to the normalized Lebesgue measure $d_Nx$ on ${\mathbb R}$
by $|\| \cdot \||_u$.\\
\ \\
{\bf Step 4.} \ {\em Apply Corollary \ref{corlieb} of Lieb theorem to compute the operatorial norm
of the bilinear operator $\tilde{B}$, by computing the supremum only over the set of exponential
functions.}\\
\ \\
We are working now in dimension $d = 1$.
Since the function $(x$, $y) \mapsto -(1/2)(mx - ny)^2$ is non--positive,
we may use Lieb theorem. That means:
\begin{eqnarray*}
& \ & \|\tilde{B}\|_{L^p({\mathbb R}, d_Nx) \times L^q({\mathbb R}, d_Nx) \to L^r({\mathbb R}, d_Nx)}\\
& = & \sup\{|\|\tilde{B}(f, g)\||_r \mid f \in L^p({\mathbb R}, d_Nx), g \in L^q({\mathbb R}, d_Nx),
|\|f\||_p = |\|g\||_q = 1\}\\
& = & \sup\{|\|\tilde{B}(f, g)\||_r \mid f = c_1\exp[-(k/2)x^2], g = c_2\exp[-(l/2)x^2],\\
& \ & \ \ \ \ \ \ \ \ \ \ \ \ \ \ \ \ \ \ \ \ \ \ \ \ c_1 > 0, c_2 > 0, k > 0, l > 0, |\|f\||_p = |\|g\||_q = 1\}.
\end{eqnarray*}
If $h(x) = \lambda\exp[-(s/2)x^2]$, with $s > 0$,
then for all $u \geq 1$, we have:
\begin{eqnarray*}
\||h|\|_u & = & |\lambda|\left[\int_{{\mathbb R}}e^{-\frac{us}{2}x^2}d_Nx\right]^{1/u}\\
({\rm let} \ x' := \sqrt{us}x) \quad & = &
|\lambda|\left[\frac{1}{\sqrt{us}}\int_{{\mathbb
R}}e^{-\frac{x'^2}{2}}d_Nx'\right]^{1/u}\\
& = & |\lambda|\frac{1}{(\sqrt{us})^{1/u}}.
\end{eqnarray*}
Thus, in order to have $\||h|\|_u = 1$, we must have:
\begin{eqnarray}
|\lambda| & = & (\sqrt{us})^{1/u}. \label{c1}
\end{eqnarray}
Therefore, we have:
\begin{eqnarray*}
& \ & \|\tilde{B}\|_{L^p({\mathbb R}, d_Nx) \times L^q({\mathbb R}, d_Nx) \to L^r({\mathbb R}, d_Nx)}\\
& = & \sup_{k > 0, l > 0}\left\{(pk)^{1/(2p)}(ql)^{1/(2q)}\right. \\
& \ & \left. \ \ \ \ \ \ \ \ \ \times \left\|\left|\int_{{\mathbb R}}f(\alpha\cdot +
\gamma y)g(\beta \cdot + \delta y)e^{-(1/2)(m\cdot - ny)^2}d_Ny\right|\right\|_r\right\} \\
& = & p^{1/(2p)}q^{1/(2q)}\sup_{k > 0, l > 0}\left\{k^{1/(2p)}l^{1/(2q)}\right. \\
& \ & \left. \ \ \ \ \times \left\{\int_{{\mathbb R}}\left[\int_{\mathbb
R}e^{-\frac{k}{2}(\alpha x + \gamma y)^2}e^{-\frac{l}{2}(\beta x + \delta y)^2}
e^{-(1/2)(mx - ny)^2}d_Ny\right]^rd_Nx\right\}^{1/r}\right\}.
\end{eqnarray*}
We apply now formula (\ref{Lagrange}) from Lemma \ref{computational_lemma} and conclude that:
\begin{eqnarray}
& \ & \int_{\mathbb
R}e^{-\frac{k}{2}(\alpha x + \gamma y)^2}e^{-\frac{l}{2}(\beta x + \delta y)^2}
e^{-(1/2)(mx - ny)^2}d_Ny \label{interior_integral_computation}\\
& = & \frac{1}{\sqrt{\gamma^2k + \delta^2l + n^2}} \nonumber\\
& \ & \times \exp\left[-\frac{1}{2} \cdot \frac{(\alpha\delta - \beta\gamma)^2kl +
(n\alpha + m\gamma)^2k + (n\beta + m\delta)^2l}{\gamma^2k + \delta^2l + n^2}x^2\right].\nonumber
\end{eqnarray}
Thus, we have:
\begin{eqnarray*}
& \ & \|\tilde{B}\|_{L^p({\mathbb R}, d_Nx) \times L^q({\mathbb R}, d_Nx) \to L^r({\mathbb R}, d_Nx)}\\
& = & p^{1/(2p)}q^{1/(2q)}\sup_{k > 0, l > 0}\left\{k^{1/(2p)}l^{1/(2q)}\right. \\
& \ & \left. \ \ \ \ \times \left\{\int_{{\mathbb R}}\left[\int_{\mathbb
R}e^{-\frac{k}{2}(\alpha x + \gamma y)^2}e^{-\frac{l}{2}(\beta x + \delta y)^2}
e^{-(1/2)(mx - ny)^2}d_Ny\right]^rd_Nx\right\}^{1/r}\right\}\\
& = & p^{1/(2p)}q^{1/(2q)}\sup_{k > 0, l > 0}\left\{k^{1/(2p)}l^{1/(2q)}
\frac{1}{\sqrt{\gamma^2k + \delta^2l + n^2}}\right.\\
& \ & \ \ \times \left\{\int_{{\mathbb R}}\exp\left[-\frac{r}{2} \cdot \frac{(\alpha\delta - \beta\gamma)^2kl +
(n\alpha + m\gamma)^2k + (n\beta + m\delta)^2l}{\gamma^2k + \delta^2l + n^2}x^2\right]\right.\\
& \ & \left.\left. \ \ \ \ \ \ \ \ \ \ \ \ \ \ \ \ \ \ \ \ \ \ \ \ \ \ \ \ \ \ \ \ \ \ \ \ \ \ \ \ d_Nx
\right\}^{1/r}\right\}.\\
& = & \frac{p^{1/(2p)}q^{1/(2q)}}{r^{1/(2r)}}\sup_{k > 0, l > 0}\left\{k^{1/(2p)}l^{1/(2q)}
\frac{1}{\sqrt{\gamma^2k + \delta^2l + n^2}}\right.\\
& \ & \left. \ \ \ \ \ \ \ \ \ \ \ \times \frac{\left(\gamma^2k + \delta^2l + n^2\right)^{1/(2r)}}
{\left[(\alpha\delta - \beta\gamma)^2kl
+ (n\alpha + m\gamma)^2k + (n\beta + m\delta)^2l\right]^{1/(2r)}}\right\}\\
& = &  \frac{p^{1/(2p)}q^{1/(2q)}}{r^{1/(2r)}}\sup_{k > 0, l > 0}\frac{k^{1/(2p)}l^{1/(2q)}}
{\left(\gamma^2k + \delta^2l + n^21\right)^{1/(2r')}\left[U^2kl
+ V^2k + W^2l\right]^{1/(2r)}},
\end{eqnarray*}
where $r'$ is the H\"older conjugate of $r$, that means:
\begin{eqnarray}
\frac{1}{r} + \frac{1}{r'} & = & 1,
\end{eqnarray}
\begin{eqnarray}
U & := & \alpha\delta - \beta\gamma,
\end{eqnarray}
\begin{eqnarray}
V & := & n\alpha + m\gamma,
\end{eqnarray}
and
\begin{eqnarray}
W & := & n\beta + m\delta.
\end{eqnarray}
We can put the factors $r$ and $U^2kl + V^2k + W^2l$ together. Thus it remains to prove that:
\begin{eqnarray}
& \ & \sup_{k > 0, l > 0}\frac{(p{\bf k})^{1/(2p)}(q{\bf l})^{1/(2q)}}
{\left(\gamma^2{\bf k} + \delta^2{\bf l} + n^2{\bf 1}\right)^{1/(2r')}
\left[rU^2{\bf kl} + rV^2{\bf k} + rW^2{\bf l}\right]^{1/(2r)}} \nonumber\\
& = & 1. \label{needed_supremum}
\end{eqnarray}
We used some bold face letters to emphasize the idea for the next step.\\
\ \\
{\bf Step 5.} \ {\em Observe that the numerator of the left--hand side of (\ref{needed_supremum})
looks like a weighted geometric mean of $(k$, $l)$, while the two factors from the denominator
of the left--hand side of (\ref{needed_supremum}) look like weighted arithmetic means of ($k$, $l$, $1$)
and ($kl$, $k$, $l$), respectively. Since each arithmetic mean dominates each geometric mean, our
supremum has great chances of being finite. This observation tells us that, to continue, we have to use
the reason behind the arithmetic--geometric mean inequality, which is the concavity of the
logarithmic function.}\\
\ \\
Let us make the following changes of variables:
\begin{eqnarray}
K & := & pk
\end{eqnarray}
and
\begin{eqnarray}
L & := & ql.
\end{eqnarray}
Thus, we have to prove that:
\begin{eqnarray}
& \ & \sup_{K > 0, L > 0}\frac{{\bf K}^{1/(2p)}{\bf L}^{1/(2q)}}
{\left[(\gamma^2/p){\bf K} + (\delta^2/q){\bf L} + n^2{\bf 1}\right]^{1/(2r')}} \nonumber\\
& \ & \times \frac{1}
{\left[(rU^2/(pq)){\bf KL} + (rV^2/p){\bf K} + (rW^2/q){\bf L}\right]^{1/(2r)}} \nonumber\\
& = & 1.
\end{eqnarray}
{\bf Claim 1:} \ We have:
\begin{eqnarray}
\frac{\gamma^2}{p} + \frac{\delta^2}{q} + n^2 & = & 1. \label{frist_sum_1}
\end{eqnarray}
Indeed, we have:
\begin{eqnarray*}
\frac{\gamma^2}{p} + \frac{\delta^2}{q} + n^2 & = & \frac{\left[{\rm sgn}(\alpha\beta(t - 1))\sqrt{1 - \alpha^2}\right]^2}{p} +
\frac{\left(\sqrt{1 - \beta^2}\right)^2}{q} + c\\
& = & \frac{1 - \alpha^2}{p} + \frac{1 - \beta^2}{q} + \left(1 - \frac{1 - \alpha^2}{p} - \frac{1 - \beta^2}{q}\right)\\
& = & 1.
\end{eqnarray*}
{\bf Claim 2.} \ We have:\\
\begin{eqnarray}
\frac{rU^2}{pq} + \frac{rV^2}{p} + \frac{rW^2}{q} & = & 1. \label{second_sum_1}
\end{eqnarray}
Indeed, we have:
\begin{eqnarray*}
& \ & \frac{rU^2}{pq} + \frac{rV^2}{p} + \frac{rW^2}{q}\\
& = & \frac{r(\alpha\delta - \beta\gamma)^2}{pq} + \frac{r(n\alpha + m\gamma)^2}{p} +
\frac{r(n\beta + m\delta)^2}{q}\\
& = & r\left[\frac{(\alpha\delta - \beta\gamma)^2}{pq} +
n^2\left(\frac{\alpha^2}{p} + \frac{\beta^2}{q}\right)
+ m^2\left(\frac{\gamma^2}{p} + \frac{\delta^2}{q}\right)
+ 2mn\left(\frac{\alpha\gamma}{p} + \frac{\beta\delta}{q}\right)\right]\\
& = & r\left[\frac{(\alpha\delta - \beta\gamma)^2}{pq} +
c\left(\frac{1}{r} - a\right)
+ a(1 - c)
+ 2\sqrt{ac} \ {\rm sgn}(b) \cdot b\right]\\
& = &  r\left[\frac{(\alpha\delta - \beta\gamma)^2}{pq} +
c\left(\frac{1}{r} - a\right)
+ a(1 - c)
+ 2\sqrt{ac} \cdot \sqrt{ac}\right]\\
& = &  r\left[\frac{\left(\alpha\sqrt{1 - \beta^2} - \beta{\rm sgn}(\alpha\beta(t - 1))\sqrt{1 - \alpha^2}\right)^2}{pq}
+ \frac{c}{r} + a\right]\\
& = &  r\left[\frac{\alpha^2\left(1 - \beta^2\right) + \beta^2\left(1 - \alpha^2\right) - 2\alpha\beta
{\rm sgn}(\alpha\beta(t - 1))\sqrt{\left(1 - \alpha^2\right)\left(1 - \beta^2\right)}}{pq}\right.\\
& \ & \ \ \ \ + \left. \frac{c}{r} + a\right].\\
\end{eqnarray*}
Using now the assumption that we have equality in condition (\ref{a_b_t_condition}), we get:
\begin{eqnarray*}
{\rm sgn}(\alpha\beta(t - 1))\sqrt{\left(1 - \alpha^2\right)\left(1 - \beta^2\right)} & = & \alpha\beta(t - 1).
\end{eqnarray*}
Thus, we obtain:
\begin{eqnarray*}
& \ & \frac{rU^2}{pq} + \frac{rV^2}{p} + \frac{rW^2}{q}\\
& = &  r\left[\frac{\alpha^2\left(1 - \beta^2\right) + \beta^2\left(1 - \alpha^2\right) - 2\alpha^2\beta^2(t - 1)}{pq}
+ \frac{c}{r} + a\right]\\
& = & r\left[\frac{\alpha^2 + \beta^2 - 2\alpha^2\beta^2t}{pq}
+ \frac{c}{r} + \frac{1}{r} - \frac{\alpha^2}{p} - \frac{\beta^2}{q}\right]\\
& = & c + 1 - \frac{r}{pq}\left[\alpha^2(q - 1) + \beta^2(p - 1) + 2\alpha^2\beta^2t\right].
\end{eqnarray*}
Since we know that:
\begin{eqnarray*}
r & = & pq\frac{c}{\alpha^2(q - 1) + \beta^2(p - 1) + 2\alpha^2\beta^2t},
\end{eqnarray*}
we obtain:
\begin{eqnarray*}
& \ & \frac{rU^2}{pq} + \frac{rV^2}{p} + \frac{rW^2}{q}\\
& = & c + 1 - c\\
& = & 1.
\end{eqnarray*}
Applying Jensen inequality to the concave function $x \mapsto \ln(x)$, we get:
\begin{eqnarray*}
\ln\left(\frac{\gamma^2}{p}K + \frac{\delta^2}{q}L + n^21\right) & \geq & \frac{\gamma^2}{p}\ln K +
\frac{\delta^2}{q}\ln L + n^2\ln 1.
\end{eqnarray*}
Exponentiating both sides of the last inequality, we obtain:
\begin{eqnarray}
\frac{\gamma^2}{p}K + \frac{\delta^2}{q}L + n^21 & \geq & K^{\gamma^2/p}L^{\delta^2/q}. \label{i1}
\end{eqnarray}
Applying again Jensen inequality, we get:
\begin{eqnarray*}
\ln\left(\frac{rU^2}{pq}KL + \frac{rV^2}{p}K + \frac{rW^2}{q}L\right) & \geq & \frac{rU^2}{pq}\ln(KL)
+ \frac{rV^2}{p}\ln K + \frac{rW^2}{q}\ln L.
\end{eqnarray*}
Exponentiating both sides of the last inequality, we obtain:
\begin{eqnarray}
\frac{rU^2}{pq}KL + \frac{rV^2}{p}K + \frac{rW^2}{q}L & \geq & K^{r[U^2/(pq) + V^2/p]}L^{r[U^2/(pq) + W^2/q]}.
\label{i2}
\end{eqnarray}
Thus applying inequalities (\ref{i1}) and (\ref{i2}), for all $K$ and $L$ positive numbers, we have:
\begin{eqnarray*}
& \ & \frac{K^{1/(2p)}L^{1/(2q)}}{\left[(\gamma^2/p)K + (\delta^2/q)L + n^21\right]^{1/(2r')}}\\
& \ & \times
\frac{1}{\left[(rU^2/(pq))KL + (r/p)V^2K + (rW^2/q)L\right]^{1/(2r)}}\\
& \leq & \frac{K^{1/(2p)}L^{1/(2q)}}{K^{\gamma^2/(2pr')}L^{\delta^2/(2qr')} \cdot K^{U^2/(2pq) + V^2/(2p)}
L^{U^2/(2pq) + W^2/(2q)}}\\
& = & \left[\frac{K}{K^{\gamma^2/r' + U^2/q + V^2}}\right]^{1/(2p)}
\left[\frac{L}{L^{\delta^2/r' + U^2/p + W^2}}\right]^{1/(2q)}.
\end{eqnarray*}
If we can prove now the following claim, then we will be done:\\
{\bf Claim 3:} \ We have:
\begin{eqnarray}
\frac{\gamma^2}{r'} + \frac{U^2}{q} + V^2 & = & 1
\end{eqnarray}
and similarly,
\begin{eqnarray}
\frac{\delta^2}{r'} + \frac{U^2}{p} + W^2 & = & 1.
\end{eqnarray}
Indeed, we have:
\begin{eqnarray*}
& \ & \frac{\gamma^2}{r'} + \frac{U^2}{q} + V^2\\
& = & \frac{\gamma^2}{r'} + \frac{(\alpha\delta - \beta\gamma)^2}{q} +
(n\alpha + m\gamma)^2\\
& = & \frac{\gamma^2}{r'} + \frac{(\alpha\delta - \beta\gamma)^2}{q} +
c\alpha^2 + 2\sqrt{ac} \ {\rm sgn}(b) \alpha\gamma + a\gamma^2\\
& = & \frac{\gamma^2}{r'} + \frac{(\alpha\delta - \beta\gamma)^2}{q}
+ \left(1 - \frac{\gamma^2}{p} - \frac{\delta^2}{q}\right)\alpha^2
+ 2b\alpha\gamma + \left(\frac{1}{r} - \frac{\alpha^2}{p} - \frac{\beta^2}{q}\right)\gamma^2\\
& = & \left(\frac{\gamma^2}{r'} + \frac{\gamma^2}{r}\right) +
\left[\frac{(\alpha\delta - \beta\gamma)^2}{q} - \frac{\alpha^2\delta^2 + \beta^2\gamma^2}{q}\right]
+ \alpha^2 - 2\frac{\alpha^2\gamma^2}{p} + 2b\alpha\gamma.
\end{eqnarray*}
We can replace now $b$ by $\alpha\gamma/p + \beta\delta/q$, and obtain:
\begin{eqnarray*}
& \ & \frac{\gamma^2}{r'} + \frac{U^2}{q} + V^2\\
& = & \gamma^2\left(\frac{1}{r'} + \frac{1}{r}\right) -
2\frac{\alpha\beta\gamma\delta}{q}
+ \alpha^2 - 2\frac{\alpha^2\gamma^2}{p} + 2\left(\frac{\alpha\gamma}{p} +
\frac{\beta\delta}{q}\right)\alpha\gamma\\
& = & \gamma^2 \cdot 1 + \alpha^2\\
& = & \left(1 - \alpha^2\right) + \alpha^2\\
& = & 1.
\end{eqnarray*}
Thus, we have proved that for all $K$ and $L$ positive, we have:
\begin{eqnarray}
& \ & \frac{K^{1/(2p)}L^{1/(2q)}}{\left[(\gamma^2/p)K + (\delta^2/q)L + C^21\right]^{1/(2r')}} \nonumber\\
& \ & \times
\frac{1}{\left[(rU^2/(pq))KL + (r/p)V^2K + (rW^2/q)L\right]^{1/(2r)}} \nonumber\\
& \leq & 1. \label{less_than_or_equal_to_1}
\end{eqnarray}
On the other hand, for $K = L = 1$, we have equality in (\ref{less_than_or_equal_to_1}). \\
Therefore, our supremum is equal to $1$.\\
\ \\
{\bf Step 6.} \ {\em We extend the inequality to the infinite dimensional case.}\\
\ \\
Let $\{e_n\}_{n \geq 0}$ be an orthonormal basis of $E$ made up of eigenfunctions of the operator $A$ used
in the construction of the Gel'fand triple ${\mathcal E} \subset E \subset {\mathcal E}'$.
For all natural numbers $n$, let:
\begin{eqnarray*}
E_n & := & {\mathbb R}\langle \cdot, e_0 \rangle \oplus {\mathbb R}\langle \cdot, e_1 \rangle \oplus \cdots \oplus
{\mathbb R}\langle \cdot, e_{n - 1} \rangle,
\end{eqnarray*}
where $\langle \cdot$, $e_i \rangle$ represents the $L^2$--normally distributed random variable generated by $e_i$,
for all $0 \leq i \leq n - 1$.
We define ${\mathcal F}_n := {\mathcal F}(E_n)$, that means,
${\mathcal F}_n$ is the smallest
sigma--algebra with respect to which $\langle \cdot$, $e_0 \rangle$, $\langle \cdot$, $e_1 \rangle$, $\dots$,
$\langle \cdot$, $e_{n - 1} \rangle$ are measurable.\\
For all $\varphi \in L^{p}({\mathcal E}'$, $\mu)$ and $n \geq 1$, we
define $\varphi_n := E[\varphi \mid {\mathcal F}_n]$, the conditional expectation of $\varphi$
with respect to ${\mathcal F}_n$.
Since $\{{\mathcal F}_n\}_{n \geq 1}$ is an increasing family of
sigma--algebras, and the sigma--algebra generated by them is
the Borel sigma--algebra ${\mathcal F}$ of ${\mathcal E}'$, from the Martingale Convergence Theorem we
have:
\begin{eqnarray*}
E[\varphi \mid {\mathcal F}_n] & \to & \varphi,
\end{eqnarray*}
as $n \to \infty$,
both almost surely and in $L^p({\mathcal E}'$, $\mu)$, for all $p \geq 1$.\\
Thus, the vector space:
\begin{eqnarray}
L_f^p({\mathcal E}') & := & \cup_{n = 1}^{\infty}L^p({\mathcal E}',{\mathcal F}_n, P)
\end{eqnarray}
is dense in $L^p({\mathcal E}'$, ${\mathcal F}$, $P)$.\\
Due to our proof in the finite dimensional case, we know that the bilinear operator
$B : L_f^p({\mathcal E}') \times L_f^q({\mathcal E}') \to L^r({\mathcal E}'$, $\mu)$,
defined by:
\begin{eqnarray}
B(\varphi, \psi) & = & \Gamma(C)\varphi \diamond_T \Gamma(D)\psi,
\end{eqnarray}
is bounded of operatorial norm $1$, and since $L_f^p({\mathcal E}') \times L_f^q({\mathcal E}')$
is dense in $L^p({\mathcal E}'$, $\mu) \times L^q({\mathcal E}'$, $\mu)$,
it has a unique bounded bilinear extension
from $L^p({\mathcal E}'$, $\mu) \times L^q({\mathcal E}'$, $\mu)$ to $L^r({\mathcal E}'$, $\mu)$.\\
\ \\
($\Leftarrow$) Let us assume now that the operator
\begin{eqnarray*}
B : L^p\left({\mathcal E}', \mu\right) \times  L^q\left({\mathcal E}', \mu\right) & \to &  L^r\left({\mathcal E}', \mu\right)
\end{eqnarray*}
defined as:
\begin{eqnarray}
B(\varphi, \psi) & = & \Gamma(C)\varphi \diamond_T \Gamma(D)\psi,
\end{eqnarray}
is bounded. Then there is a constant $k > 0$, such that, for all $\varphi \in L^p({\mathcal E}'$, $\mu)$ and all
$\psi \in L^q({\mathcal E}'$, $\mu)$, we have:
\begin{eqnarray}
\parallel \Gamma(C)\varphi \diamond_T \Gamma(D)\psi \parallel_r & \leq & k\parallel \varphi \parallel_p\parallel \psi \parallel_q.
\label{k_inequality}
\end{eqnarray}
Since $C$, $D$, and $T$ are diagonalized in the same base, $\{e_i\}_{i \geq 0}$ as $A$, for each $i \geq 0$, there exists $\alpha_i$,
$\beta_i$, and $t_i$ real numbers such that:
\begin{eqnarray}
Ce_i & = & \alpha_ie_i,
\end{eqnarray}
\begin{eqnarray}
De_i & = & \beta_ie_i,
\end{eqnarray}
\begin{eqnarray}
Te_i & = & t_ie_i.
\end{eqnarray}
Let $i \geq 0$ be a fixed natural number. Let $u$ and $s$ be arbitrary real numbers,
such that $u \neq 0$.
Let $\varphi$ and $\psi$ be the following exponential functions:
\begin{eqnarray}
\varphi & := & \varphi_{sue_i}
\end{eqnarray}
and
\begin{eqnarray}
\psi & := & \varphi_{ue_i}.
\end{eqnarray}
Then inequality (\ref{k_inequality}) becomes:
\begin{eqnarray}
\parallel e^{t_i\alpha_i\beta_isu^2}\varphi_{(s\alpha + \beta)ue_i}\parallel_{r} & \leq & k\parallel \varphi_{sue_i} \parallel_{p}
\parallel \varphi_{ue_i} \parallel_{q}. \label{su_inequality}
\end{eqnarray}
Since a simple computation shows that for every $l \in [1$, $\infty)$, and every exponential function $\varphi_{\xi}$, with
$\xi \in E$, we have:
\begin{eqnarray}
\parallel \varphi_{\xi} \parallel_l & = & e^{(l - 1)|\xi|_0^2/2},
\end{eqnarray}
inequality (\ref{su_inequality}) becomes:
\begin{eqnarray}
& \ & \exp\left(\frac{1}{2}(r - 1)(s\alpha_i + \beta_i)^2u^2 + t_i\alpha_i\beta_iu^2\right) \nonumber\\
& \leq & k\exp\left(\frac{1}{2}(p - 1)s^2u^2 + \frac{1}{2}(q - 1)u^2\right).
\end{eqnarray}
Taking first $\ln$ from both sides of the last inequality, and then dividing both sides by $u^2$, we obtain:
\begin{eqnarray}
\frac{1}{2}(r - 1)(s\alpha_i + \beta_i)^2 + t_i\alpha_i\beta_i
& \leq & \frac{\ln k}{u^2} + \frac{1}{2}(p - 1)s^2 + \frac{1}{2}(q - 1),
\end{eqnarray}
for all $u \neq 0$ and $s \in {\mathbb R}$. Passing to the limit, as $u \to \infty$, in this inequality, we obtain:
\begin{eqnarray}
\frac{1}{2}(r - 1)(s\alpha_i + \beta_i)^2 + t_i\alpha_i\beta_i
& \leq & \frac{1}{2}(p - 1)s^2 + \frac{1}{2}(q - 1), \label{s_inequality}
\end{eqnarray}
for all real numbers $s$. Inequality (\ref{s_inequality}) is equivalent to:
\begin{eqnarray*}
\frac{1}{2}\left[p - 1 - \alpha_i^2(r - 1)\right]s^2 - \left(t_i + r - 1\right)\alpha_i\beta_is + \frac{1}{2}\left[q - 1 - \beta_i^2(r - 1)\right]
& \geq & 0,
\end{eqnarray*}
for all real numbers $s$. For this quadratic function of variable $s$ to be non--negative, for all real values of $s$, its leading coefficient
$(1/2)[p - 1 - \alpha_i^2(r - 1)]$ must be non--negative and its discriminant must be non--positive, that means:
\begin{eqnarray*}
\alpha_i^2\beta_i^2\left(t_i + r - 1\right)^2 - \left[p - 1 - \alpha_i^2(r - 1)\right]\left[q - 1 - \beta_i^2(r - 1)\right] & \leq & 0.
\end{eqnarray*}
The last inequality is equivalent to:
\begin{eqnarray}
r - 1 & \leq & \frac{(p - 1)(q - 1) - t_i^2\alpha_i^2\beta_i^2}{\alpha_i^2(q - 1) + \beta_i^2(p - 1) + 2\alpha_i^2\beta_i^2t_i}.
\end{eqnarray}
Since this inequality holds for all $i \geq 0$, we conclude that:
\begin{eqnarray}
(r - 1)I & \leq & \frac{(p - 1)(q - 1)I - C^2D^2T^2}{(q - 1)C^2 + (p - 1)D^2 + 2C^2D^2T}.
\end{eqnarray}
\end{proof}

\noindent Let us define for any invertible self--adjoint operator $B$ on $E_c$ commuting with $A$, and any $p \in [1$, $\infty]$, the following norm:
\begin{eqnarray}
\parallel \varphi \parallel_{p, B} & := & \parallel \Gamma(B)\varphi \parallel_p.
\end{eqnarray}
Let us also define the space:
\begin{eqnarray}
L^{p, B}({\mathcal E}', \mu) & := & \{\varphi \in ({\mathcal E})^* ~\mid~ \parallel \Gamma(B)\varphi \parallel_p ~<~ \infty\}.
\end{eqnarray}
With this notation, we have the following corollary.
\begin{corollary} \label{true_Holder_theorem}
Let ${\mathcal E} \subset E \subset {\mathcal E}'$ be a Gel'fand triple
given by a self--adjoint diagonal operator $A$ on $E$, with increasing, greater
than $1$ eigenvalues, whose inverse is a Hilbert--Schmidt operator. Let
$\mu$ be the Gaussian probability measure on ${\mathcal E}'$ whose existence
is guaranteed by Minlos Theorem.
Let $T$ be a self--adjoint, diagonal operator on $E$, commuting with $A$, such that:
\begin{eqnarray}
T & \geq & 0.
\end{eqnarray}
Let $B$, $C$, and $D$ be three invertible self--adjoint and diagonal operators
on $E$, commuting with the operator $A$, such that:
\begin{eqnarray}
|B| & \geq & \sqrt{\frac{T}{2}}, \label{condition_B_T}
\end{eqnarray}
\begin{eqnarray}
|C| & \geq & |B|, \label{condition_C_B}
\end{eqnarray}
\begin{eqnarray}
|D| & \geq & |B|, \label{condition_D_B}
\end{eqnarray}
and
\begin{eqnarray}
\left(C^2 - B^2\right)(D^2 - B^2) & \geq & \left(T - B^2\right)^2. \label{B_C_D_T_condition}
\end{eqnarray}
Let $p$, $q$, $r > 1$ such that:\\
\begin{eqnarray}
\frac{1}{(r - 1)B^2 + T} & \geq & \frac{1}{(p - 1)C^2 + T} +
\frac{1}{(q - 1)D^2 + T}, \label{true_Holder_condition}
\end{eqnarray}
Then for all $\varphi$ in $L^{p, C}({\mathcal E}'$, $\mu)$ and
$\psi$ in $L^{q, D}({\mathcal E}'$, $\mu)$,
$\varphi \diamond_T \psi$ belongs
to $L^{r, B}({\mathcal E}'$, $\mu)$, and the following
inequality holds:
\begin{eqnarray}
\parallel \varphi \diamond_T \psi \parallel_{r, B} & \leq & \parallel \varphi \parallel_{p, C}
\cdot \parallel \psi \parallel_{q, D}. \label{true_Holder_inequality}
\end{eqnarray}
On the other hand, if:
\begin{eqnarray*}
\frac{1}{(r - 1)B^2 + T} & \not{\geq} & \frac{1}{(p - 1)C^2 + T} +
\frac{1}{(q - 1)D^2 + T},
\end{eqnarray*}
then the bilinear operator $(\varphi$, $\psi) \mapsto \varphi \diamond_T \psi$
is not bounded from  $L^{p, C}({\mathcal E}'$, $\mu) \times L^{q, D}({\mathcal E}'$, $\mu)$ to
$L^{r, B}({\mathcal E}'$, $\mu)$.
\end{corollary}
\begin{proof}
Simply apply Theorem \ref{holder_young_theorem} and Lemma \ref{functorial_lemma} to the following operators and
random variables, with the convention
that $X \rightarrow X'$ means that $X$ from theorem \ref{holder_young_theorem} is replaced by $X'$ in the same theorem:
\begin{itemize}
\item $C \rightarrow BC^{-1}$

\item $D \rightarrow BD^{-1}$

\item $T \rightarrow TB^{-2}$

\item $\varphi \rightarrow \Gamma(C)\varphi$

\item $\psi \rightarrow \Gamma(D)\psi$
\end{itemize}
Writing inequality (\ref{fullholder}) for these new operators and random variables, and using the
fact that:
\begin{eqnarray}
\parallel \varphi \parallel_{p, |B|} & = & \parallel \varphi \parallel_{p, B},
\end{eqnarray}
inequality (\ref{true_Holder_inequality}) follows.
\end{proof}
\begin{corollary}
If we choose $B := I$, $C > I$, $\psi := 1$, and let either the eigenvalues of $D$ go to $\infty$, or choose the
eigenvalues of $D$ large enough such that
(\ref{B_C_D_T_condition}) is satisfied and let $q$ go to $\infty$, then condition (\ref{true_Holder_condition})
becomes:
\begin{eqnarray}
\frac{1}{(r - 1)I + T} & \geq & \frac{1}{(p - 1)C^2 + T} + \frac{1}{\infty},
\end{eqnarray}
which is equivalent to:
\begin{eqnarray}
|C| & \geq & \sqrt{\frac{r - 1}{p - 1}}I \label{Nelson_condition}
\end{eqnarray}
and inequality (\ref{true_Holder_inequality}) becomes:
\begin{eqnarray}
\parallel \varphi \parallel_r & \leq & \parallel \Gamma(C)\varphi \parallel_p. \label{Nelson_inequality}
\end{eqnarray}
Inequalities (\ref{Nelson_condition}) and (\ref{Nelson_inequality}) are exactly Nelson condition and
hypercontractivity inequality.
\end{corollary}

\begin{corollary}
If we choose $B = I$, then condition (\ref{condition_B_T}) becomes:
\begin{eqnarray}
0 \leq T \leq 2I,
\end{eqnarray}
and conditions (\ref{condition_C_B}) and (\ref{condition_D_B}) become:
\begin{eqnarray}
|C| & \geq & I
\end{eqnarray}
and
\begin{eqnarray}
|D| & \geq & I.
\end{eqnarray}
In this case the inequality:
\begin{eqnarray}
\frac{1}{(r - 1)I + T} & \geq & \frac{1}{(p - 1)C^2 + T} +
\frac{1}{(q - 1)D^2 + T},
\end{eqnarray}
and the ``smoothness" conditions $\Gamma(C)\varphi \in L^p({\mathcal E}'$, $\mu)$
and $\Gamma(D)\varphi \in L^q({\mathcal E}'$, $\mu)$ guarantee the fact that
$\varphi \diamond_T \psi$ is a \underline{true random variable}
(not a merely generalized function)
in the space $L^r({\mathcal E}'$, $\mu)$.
\end{corollary}

We would like to make the following comments:\\
\par {\bf Comments:}

\begin{itemize}

\item Inequality (\ref{fullholder}) in the case $T := 0$, $C := \alpha I$, $D := \beta I$, $\alpha^2 + \beta^2 = 1$, and $p = q = r := 2$,
was established using Cauchy--Schwarz inequality in \cite{kss2002}.

\item Inequality (\ref{fullholder}) in the case $T := 0$, $C := \alpha I$, $D := \beta I$, $\alpha^2 + \beta^2 = 1$, and
($p = q = r := 1$ or $p = q = r := \infty$) was proven in \cite{ls2008}.

\item Inequality (\ref{fullholder}) in the case $T := 0$, $C := \alpha I$, $D := \beta I$, $\alpha^2 + \beta^2 = 1$, and
$1/(r - 1) = \alpha^2/(p - 1) + \beta^2/(q - 1)$ was proven in \cite{dpls11}.

\item In the case: $B = C = D = T = I$, the $T$--Wick product becomes the point--wise product, and
condition (\ref{true_Holder_condition}) becomes:
\begin{eqnarray}
\frac{1}{r} & \geq & \frac{1}{p} + \frac{1}{q},
\end{eqnarray}
which is exactly the classic H\"older condition for probability measures. One should not forget, that
for general measures, H\"older condition is the equality:
\begin{eqnarray}
\frac{1}{r_0} & = & \frac{1}{p} + \frac{1}{q},
\end{eqnarray}
not an inequality, but if the measure is a probability measure, Lyapunov inequality:
for all $0 < r \leq r_0$, we have:
\begin{eqnarray}
\parallel f \parallel_r & \leq & \parallel f \parallel_{r_0},
\end{eqnarray}
relaxes the H\"older condition from a perfect equality to an inequality.

\item Methods of finding the supremum over the exponential functions in Lieb theorem were provided in
\cite{bcct05}, \cite{cce09}, and \cite{cll04}. In this paper we used a method based on Jensen inequality
for the concave natural logarithmic function that was also applied in \cite{dpls11}.

\end{itemize}

\end{document}